\documentclass[twoside]{siamltex}
\usepackage{amssymb}
\usepackage{amsmath}
\usepackage[final]{graphicx}
\usepackage{psfrag}
\usepackage{epsfig}
\usepackage{latexsym}
\usepackage[final]{showkeys}
\usepackage{exscale}
\usepackage{algorithm}
\usepackage{algorithmic}

\font\Bbb=msbm10 

\newcommand{\Rset}{\mbox{\Bbb R}}

\newcommand{\weg}[1]{ }
\newcommand{\Field}[1]{\mbox{\boldmath$#1$\unboldmath}}
\newcommand{\Tensor}[1]{\mbox{\boldmath$#1$\unboldmath}}
\newcommand{\Svec}[1]{\mbox{\boldmath$#1$\unboldmath}}

\def \veps{\varepsilon}
\def \vfi{\varphi}
\def \eps{\epsilon}
\def \half{{\textstyle\frac12}}
\title{Adaptive, Fast and Oblivious Convolution \\
in Evolution Equations with Memory}
\author{Mar\'ia L\'opez-Fern\'andez\footnotemark[1] \and
Christian Lubich\footnotemark[5] \and
Achim Sch\"adle\footnotemark[7]}
\begin{document}
\maketitle
\renewcommand{\thefootnote}{\fnsymbol{footnote}}
\footnotetext[1]{Departamento de Matem\'atica Aplicada,
Universidad de Valladolid, Valladolid, Spain.
 ~E-mail: {\tt marial@mac.cie.uva.es}. Supported by DGI-MCYT under
 project MTM 2004-07194 cofinanced by FEDER funds.}
\footnotetext[5]{Mathematisches Institut,
  Universit\"at T\"ubingen,
  Auf der Morgenstelle 10,
  D--72076 T\"ubingen,
  Germany.
  ~
  E-mail:  {\tt lubich@na.uni-tuebingen.de}. Supported by DFG, SFB 382.}
\footnotetext[7]{ZIB Berlin, Takustr.~7, D-14195 Berlin, Germany.
~E-mail: {\tt schaedle@zib.de}. Supported by the DFG Research
Center \textsc{Matheon} ``Mathematics for key technologies'', Berlin.}
\renewcommand{\thefootnote}{\arabic{footnote}}
\begin{abstract}
To approximate convolutions which occur in
evolution equations with memory terms, a variable-stepsize algorithm is presented for which
advancing $N$ steps requires only $O(N\log N)$ operations and
$O(\log N)$ active memory, in place of $O(N^2)$ operations
and $O(N)$ memory for a direct implementation.
A basic feature of the fast algorithm is the reduction, via contour integral
representations, to differential equations  which are solved
numerically with adaptive step sizes.
Rather than the kernel itself, its Laplace transform is used
in the algorithm.
The algorithm is illustrated on three examples:
a blow-up example originating from a Schr\"odinger
equation with concentrated nonlinearity, chemical reactions with inhibited diffusion, and
viscoelasticity with a fractional order constitutive law.
\end{abstract}
\section{Introduction}
\label{Sec:Intro}
We consider the problems of computing the convolution
\begin{equation}\label{eq:convolution}
\int_0^t f(t-\tau)\, g(\tau)\, d\tau\,, \qquad 0\le t\le T,
\end{equation}
possibly with matrix-valued kernel $f$ and vector-valued function
$g$,
and of solving evolution equations with memory containing such
convolution integrals where
$g$ is not a function known in advance, but $g(\tau)$ depends
on the solution at time $\tau$ of the integral equation or
integro-differential equation.
In previous papers \cite{LuS02,SchLoLu06}
we have developed convolution algorithms that are
{\it fast} and {\it oblivious}: to approximate (\ref{eq:convolution}) on a
grid $t=nh$ $(n=0,1,\dots,N)$ with constant step size $h$ and $Nh=T$,
the algorithm requires
\begin{itemize}
\item $O(N\log N)$ operations,
\item $O(\log N)$ evaluations
of the Laplace transform $F={\mathcal L}f$, none of $f$, and
\item $O(\log N)$
active memory.
\end{itemize}
In the $n$th time step, $g$ is evaluated at $t_n=nh$, but the
history $g(t_j)$ for $j<n$ is forgotten in this algorithm, and only
logarithmically few linear combinations of the values of $g$ are
kept in memory. This is to be contrasted with the $O(N^2)$
operations, $O(N)$ evaluations of the kernel $f$, and $O(N)$ memory
for a naive implementation of a quadrature formula for
(\ref{eq:convolution}). Moreover, we note that in many applications
the Laplace transform $F$ (the transfer function), rather than the
kernel $f$ (the impulse response), is known {\it a priori}. A basic
feature of the fast algorithm is the reduction, via contour integral
representations, to differential equations of the form $y'=\lambda y
+ g$ for suitable complex values of $\lambda$, which are solved
numerically. It is not  necessary to solve these differential
equations with constant time step $h$, as was done in
\cite{LuS02,SchLoLu06}, but the step size may instead be adapted to
the behavior of $g$. This observation opens the way to an {\it
adaptive} fast and oblivious convolution algorithm. Turning this
simple idea into an efficient algorithm is, however, not as simple
and the development of such an algorithm is precisely the topic of
the present paper. The need to use adaptive time steps in solving
evolutionary integro-differential equations in applications has been
addressed at various places in the literature, e.g., by Adolfsson,
Enelund \& Larsson \cite{AdoEL04}, Cao, Burrage \& Abdullah
\cite{CaoBA06}, and Diethelm \& Freed~\cite{DF98}. None of the
adaptive algorithms proposed there, however, can make use of the
convolution structure to reduce the $O(N^2)$ operation count and
$O(N)$ memory requirements for $N$ steps. The convolution algorithm
proposed here works in the situation of a {\it sectorial} Laplace
transform $F$:
\begin{equation}
  \label{eq:sector}
  \begin{array}{c}
    \hbox{$F$ is analytic in a sector $|\arg(s-\sigma)|<\pi-\varphi$ with
      $\varphi<\half\pi$, 
      }
    \\[1mm]
  \hbox{and in this sector, }\
    |F(s)| \le M\, |s|^{-\nu}\quad\hbox{for some real $M$ and $\nu>0$}.
  \end{array}
\end{equation}
An equivalent condition is that $f$ is analytic in a complex sector
containing the positive real half-axis $t>0$, and is bounded by
$|f(t)|\le C\, t^{\nu-1}e^{\sigma t}$ in this sector.
A typical example is the fractional-power kernel
$f(t)=t^{\nu-1}/\Gamma(\nu)$, which has the Laplace
transform $F(s)=s^{-\nu}$.
An essential ingredient of the algorithm is the discretization of
the inversion formula for the Laplace transform, given by
\begin{equation}\label{eq:inv}
  f(t) = \frac{1}{2\pi i} \int_\Gamma   e^{t\lambda }\,F(\lambda)\,  d\lambda,
  \qquad t>0,
\end{equation}
with $\Gamma$ a contour in the sector of analyticity oriented with
increasing imaginary part and going to infinity
with an acute angle to the negative real half-axis, so that
$e^{t\lambda}$ decays fast for growing $|\lambda|$ along $\Gamma$.
We will choose the contour as a hyperbola. Since we
cannot obtain a uniformly good approximation for all $t\in(0,T]$
with a single contour~$\Gamma$,
we use different hyperbolas
$\Gamma_\ell$ corresponding to geometrically growing
intervals of uniform approximation, $t\in [B^{\ell-1}h_*,B^{\ell+1}h_*]$
with an integer $B\ge 1$, e.g., $B=5$, and with a minimum step size $h_*$.
The required number of contours is thus bounded by $\log_B (T/h_*)$.
This logarithm shows up in the complexity estimates in place of
$\log_B(N)$ for the fixed-stepsize algorithm. As in that case, it appears
multiplied with the number of quadrature points on each hyperbola,
which is $O(\log\frac1\veps)$  for an accuracy $\veps$ in
the approximation of (\ref{eq:inv}).
In \S 2 we briefly review
recent results from \cite{LoP04,LoPSch} on the approximation of
inverse Laplace transforms by discretized contour integrals.
In \S 3 we describe the fast and oblivious
algorithm for computing convolutions with
variable time steps.
The algorithm is then illustrated on various problems where adaptive
time steps are important:  a blow-up
problem for a nonlinear Abel integral equation resulting from a
nonlinear Schr\"odinger equation with concentrated nonlinearity
(\S 4),  a fractional diffusion-reaction system from chemical reaction
kinetics (\S 5),
and visco-elasticity with a weakly singular memory kernel in the
constitutive equations, under applied
forces that are switched on and off~(\S 6).
\section{Preparation: Numerical inversion of the Laplace transform}
\label{Sec:InvLaplace}
In the inversion formula (\ref{eq:inv}) we choose
$\Gamma$ as the left branch of a hyperbola parameterized by
\begin{eqnarray}\label{hyppara}
  \Rset \to \Gamma &:& \ x \mapsto \gamma(x) =
  \mu (1 - \sin(a+ix)) + \sigma,
\end{eqnarray}
where $\mu > 0$ is a scale parameter, $\sigma$ is the shift in
\eqref{eq:sector}, and $0 < a < \pi/ 2 -\varphi$. Thus, $\Gamma$ is
the left branch of the hyperbola with center at $(\mu,0)$, foci at
$(0,0)$, $(2\mu,0)$, and with asymptotes forming angles $\pm(\pi/2 +
a)$ with the real axis, so that $\Gamma$ remains in the sector
\eqref{eq:sector} of analyticity of $F$. After parameterizing
\eqref{eq:inv}, the function $f$ is approximated by applying the
truncated trapezoidal rule to the resulting integral along the real
axis, i.e.,
\begin{eqnarray}\label{numinvLT}
f(t) = \frac{1}{2\pi i} \int_\Gamma   e^{t\lambda }\,F(\lambda)\,
d\lambda
\approx \sum_{k=-K}^K w_k\, e^{t\lambda_k}\,F(\lambda_k),
\end{eqnarray}
with weights $w_k$ and quadrature nodes $\lambda_k$ given by
$$
w_k = {\tau \over 2 \pi i}\: \gamma'(k\tau)~, \quad \lambda_k =
\gamma(k\tau)
$$
and $\tau > 0$ a suitable step length parameter. Different choices
of contours $\Gamma$ and parameterizations have been studied for the
numerical inversion of  sectorial Laplace transforms in the last
years. The choice of a hyperbola has been studied in
\cite{LoP04,LoPSch,Mcl,She,GavHaKh04,GavHaKh05,WeiT}, and actually we
follow here the approach in \cite{LoP04,LoPSch}. The choice of
$\Gamma$ as a parabola has also been considered recently in
\cite{GavHaKh04,GavHaKh05,WeiT}. Finally we refer to Talbot's method
\cite{Talbot79,Riz,Wei}, which could also be used with the present
algorithm; cf.~\cite{Sch02,LuS02,SchLoLu06,HiS04}. The good behavior
of this quadrature formula to approximate \eqref{eq:inv} is due to
the good properties of the trapezoidal rule when the integrand can
be analytically extended to a horizontal strip around the real axis
\cite{Ste,Sten}. We refer to \cite{LoP04,LoPSch} for details and
only give here the following error bound, which decays exponentially
in the number of quadrature nodes.
\begin{theorem}
\label{thm:err} \cite{LoPSch} Suppose that the Laplace transform $F$
satisfies the sectorial condition \eqref{eq:sector}. For fixed $T
>0$, $\Lambda \ge 1$,  $0 <a < \pi/2-\vfi$, and $K\ge 1$ there are
positive numbers
$C_1,C_2,C,c$ depending on $a$ and $\Lambda$ ($C$ depends
additionally on $T$ unless $\sigma<0$ in \eqref{eq:sector}) such
that the choice of parameters $ \tau = C_1/K$ and $ \mu = C_2
K/(\Lambda t_0)$ yields a quadrature error in \eqref{numinvLT}
bounded by
\begin{eqnarray*}
|E_K(t)| \le C\,t^{\nu-1}\, \big( \eps + e^{-cK} \big),
\end{eqnarray*}
uniformly for $t$ in intervals $[t_0,\Lambda t_0]$
with arbitrary $0<t_0\le T /\Lambda$, where $\nu$ is the exponent
of \eqref{eq:sector} and $\eps$ is the precision in the evaluations
of the Laplace transform $F$ and the elementary operations
in~\eqref{numinvLT}.
\end{theorem}
\\
Hence, $K=O(\log\frac1\veps)$ quadrature points are sufficient to
obtain an accuracy $O(\veps)$ in the approximation of the contour
integral. In practice, we choose $a \approx \half(\frac
{\pi}2-\vfi)$ and compute the values $C_1$ and $C_2$ following the
optimization process described in~\cite{LoPSch}.
\section{The variable-step-size, fast and oblivious convolution
algorithm}\label{Sec:alg}
\subsection{Local reduction to differential equations}
We want to approximate
\begin{equation}\label{eq:conv}
u(t) = \int_0^t f(t-\tau)\, g(\tau)\, d\tau
\end{equation}
on a sequence of times $0<t_1<\dots <t_N$,
spaced arbitrarily. For the moment we assume that $g$ is a known function,
though we will see later how to use the algorithm for solving
integral and integro-differential equations.
For a given $t_n$, we can insert
the Laplace inversion formula in \eqref{eq:conv} and write
\begin{equation}\label{eq:insertLTinv}
\int_0^{t_n} f(t_{n}-\tau) \, g(\tau)\,d\tau =
\int_0^{t_n} \frac 1{2\pi i}\int_{\Gamma} e^{(t_{n}-\tau)\lambda}
F(\lambda)\,d\lambda \, g(\tau)\,d\tau\,.
\end{equation}
The numerical inversion of the Laplace transform is performed very
efficiently by  the quadrature rule (\ref{numinvLT}). In
Section~\ref{Sec:InvLaplace} we have seen that the same contour
$\Gamma$ used in this quadrature can be used to approximate $f$ at
different values of $t$, ranging over intervals of the form $[t_0,
\Lambda t_0]$, for a given ratio $\Lambda \ge 1$. Since in
\eqref{eq:insertLTinv} we need to approximate $f(t_n-\tau)$ for
$t_n-\tau \in [0,t_n]$, we cannot use a unique contour $\Gamma$ and
we need to split the integral in \eqref{eq:conv} into several
pieces.
For suitable intermediate times $0 < t^- < t^+ < t_n$, with
$(t_n-t^-)/(t_n-t^+) \le \Lambda$, we select a suitable contour
$\Gamma$ for the time interval $[t_n-t^-,t_n-t^+]$ and approximate
\begin{eqnarray}\label{eq:reduction_ode}
  \int_{t^-}^{t^+} f(t_{n}-\tau) g(\tau)\,d\tau &=&
  \int_{t^-}^{t^+}  \frac 1{2\pi i}\int_{\Gamma} e^{(t_{n}-\tau)\lambda}
  F(\lambda)\, d\lambda \, g(\tau)\,d\tau \nonumber \\
  &\approx& \int_{t^-}^{t^+} \sum_{k=-K}^K w_k\, e^{(t_n-\tau)\lambda_k}\,
  F(\lambda_k)\, g(\tau)\,d\tau \nonumber \\
  &=& \sum_{k=-K}^K w_k\, F(\lambda_k) \, e^{(t_n-t^+)\lambda_k}
  \int_{t^-}^{t^+}
  e^{(t^{+}-\tau) \lambda_k}  g(\tau)\,d\tau  \nonumber \\
  &=& \sum_{k=-K}^K w_k\, F(\lambda_k)\,  e^{(t_n-t^+)\lambda_k}\,
  y(t^+,t^-,\lambda_k),
\end{eqnarray}
where $y(t^+,t^-,\lambda_k)$ is the solution at $t^{+}$ to the
linear inhomogeneous ODE
\begin{equation}\label{eq:linode}
y' = \lambda_k y + g, \quad y(t^{-}) = 0,\qquad -K\le k \le K.
\end{equation}
We now approximate $y(t^+,t^-,\lambda_k)$ by interpolating $g$
linearly on each interval $[t_{j}, t_{j+1}]$ for $j = 0,\dots,N-1$,
and integrating exactly. (More elaborate integration methods could
be used instead, cf.~\cite{ConDP06,SchLoLu06}, but for simplicity of
presentation we will just consider this particular integration
scheme.) We denote by $\bar g$ the interpolant of $g$ and by $\bar
y(t^{+},t^{-},\lambda_k)$ the resulting approximation to
$y(t^{+},t^{-},\lambda_k)$, i.e., $\bar y(t^{+},t^{-},\lambda_k)$ is
the exact solution at $t^+$ to
\begin{gather}
  \label{eq:ode}
  \begin{aligned}
    y' = \lambda_k y + \bar{g}, \quad
    y(t^{-}) = 0,\qquad -K\le k \le K.
  \end{aligned}
\end{gather}
Thus, we approximate
\begin{equation}\label{eq:approx_convode}
\int_{t^-}^{t^+} f(t_{n}-\tau) \,g(\tau)\,d\tau \approx \sum_{k=-K}^K
w_k\, F(\lambda_k)\,  e^{(t_n-t^+)\lambda_k}\, \bar y(t^+,t^-,\lambda_k)\,.
\end{equation}
\subsection{Filling the mosaic}\label{Subsec:filling}
\begin{figure}[t]
  \centering
  \includegraphics[width=0.75\textwidth]{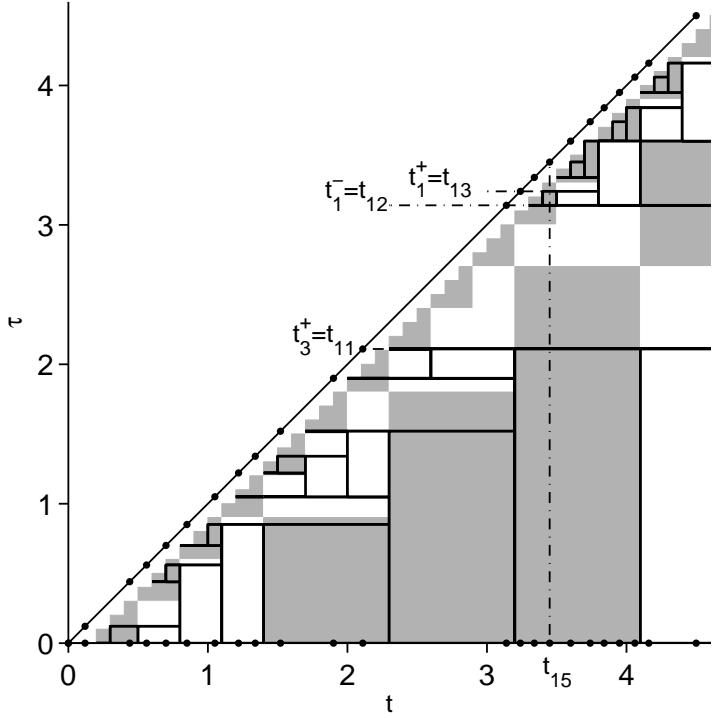}
  \caption{Mosaic in the triangle $\tau \le t$ for $B=3$ with times $t_j$
  indicated by points. Each monochromatic
  rectangle fully enclosed by black lines corresponds to
  a solution value $y(t_{\ell}^+,t_{\ell}^-,\lambda)$ of a linear differential equation
  (\ref{eq:linode}).}
  \label{fig:mosaicYM}
\end{figure}

The key to the algorithm is the way the splitting times $t^{\pm}$
for the integral in \eqref{eq:conv} are selected for every $t_n$
with $1\le n \le N$. This is done following the mosaic in the
triangle $\{(t,\tau):0\le \tau\le t \le T\}$ shown in
Fig.~\ref{fig:mosaicYM}, where patches grow geometrically with
increasing distance from the diagonal. For the moment, we fix a
minimum size of the patches closest to the diagonal, corresponding
to a minimum step size $h_*$. If along the vertical line at $t_n$
joining $0$ with the diagonal value $t_n$ we have $L$ different
patches of the tessellation, then we obtain the values $t^-_{\ell}$
and $t^+_{\ell}$ for $1 \le \ell \le L$ as the smallest and largest
grid points, resp., within the $\ell$th patch along the vertical
line at $t_n$. All those $\ell$s are collected in an index set $J$.
In case a patch does not contain any grid points its value $\ell$ is
not contained in $J$. The times $t^{\pm}_{\ell}=t^{\pm}_{n,\ell}$
depend on $n$, though for simplicity we omit this dependence in the
notation. Each class of patches of the same size in the mosaic
represents a distance class to the diagonal in the mosaic, and thus
corresponds to a different approximation interval and a different
contour to perform the inversion of the Laplace transform and to a
different set of $2K+1$ scalar differential equations. The
approximation intervals for the values $t_n-t^{\pm}_{\ell}$ are of
the form $I_{\ell}=[B^{\ell-1}h_*,B^{\ell+1}h_*]$, $2\le \ell\le L$,
so that the ratio $\Lambda$ is
$B^2$.
Since we  consider a non-equidistant sequence of times $t_j$, in
this splitting there likely appear ``gaps'' in between the
$t^+_{\ell+m}$ and $t^-_{\ell}$, which in Figure~\ref{fig:mosaicYM}
correspond to pairs of horizontal lines with any boundary of $m$
patches in between them.
For example, at the time point $t_{15}=3.45$ we have
$L=3$, $J=\{1,3\}$, $t^{-}_{3}=0$, $t^{+}_{3}=t_{11}=2.11$,
$t^{-}_{1}=t_{12}=3.14$ and $t^{+}_{1}=t_{13}=3.24$.

We split \eqref{eq:conv} into
$2|J|+1$ parts
\begin{equation}
  u(t_{n}) =
  \tilde{u}^{(0)}(t_{n}) +
  \sum_{\ell\in J} u^{(\ell)}(t_{n}) +
  \sum_{\ell\in J} \tilde{u}^{(\ell)}(t_{n})
  \label{eq:convsplit}
\end{equation}
where
\begin{equation}
  u^{(\ell)}(t_{n}) = \int_{t^{-}_{\ell}}^{t^{+}_{\ell}} f(t_{n}-\tau)
  g(\tau)\; d\tau
  \label{eq:convode}
\end{equation}
is computed using ~\eqref{eq:reduction_ode} and
\begin{equation}
  \tilde{u}^{(0)}(t_{n}) =
  \int_{t_{n-1}}^{t_{n}} f(t_{n}-\tau) g(\tau)\; d\tau \quad
  \mbox{and} \quad \tilde{u}^{(\ell)} (t_{n})=
  \int_{t^{+}_{\ell+m}}^{t^{-}_{\ell}} f(t_{n}-\tau) g(\tau)\; d\tau
  \label{eq:convds}
\end{equation}
correspond respectively to the step from $t_{n-1}$ to $t_n$ near the
diagonal and to the gaps between $t^+_{\ell+m}$ and $t^-_{\ell}$;
see Fig.~\ref{fig:mosaicYM}. These parts are computed by ``direct
steps'' explained in the next subsection. Thus, for $t_{15}=3.45$,
$u(t_{15})$ is calculated using:
\begin{itemize}
\item one  ``ode step'', $u^{(3)}$, from $t^{-}_{3}=t_{0}=0$ to  $t^{+}_{3}=t_{11}$.
\item one ``direct step'', $\tilde{u}^{(3)}$, from $t^{+}_{3}=t_{11}$ to $t^{-}_{1}=t_{12}$.
\item one ``ode step'', $u^{(1)}$, from $t^{-}_{1}=t_{12}$ to
$t^{+}_{1}=t_{13}$.
\item one ``direct step'', $\tilde{u}^{(1)}$, from $t^{+}_{1}=t_{13}$ to
$t_{14}$.
\item one ``direct step'', $\tilde{u}^{(0)}$, from $t_{14}$ to  $t_{15}$.
\end{itemize}
Translating the above splitting from the picture into a formal
procedure, we thus proceed as follows: given a minimum step size
$h_*$ and a basis $B$, for each $t_{n}$ we take $L$ as small as
possible so that we can represent $\lceil t_{n}/h_* \rceil =
2+\sum_{\ell=1}^{L} b_{\ell} B^{\ell-1}$ with $b_{\ell} \in
\{1,2,\dots,B \}$, where $\lceil x \rceil$ denotes the smallest
integer greater than $x$. $t^{+}_{\ell}$ is the largest and
$t^{-}_{\ell}$ is the smallest value in $\{t_{j}:j = 0,\dots,n\}$,
such that
$$
h_*\sum_{k=\ell+1}^{L} b_{k} B^{k-1}
\le
t_{\ell}^{-} \le t_{\ell}^{+} \le
h_*\sum_{k=\ell}^{L} b_{k} B^{k-1}.
$$
We remark that with this definition, there always exists some $j\in
\{0,\dots,n-1\}$ and some integer $m>0$ such that
$t_{\ell+m}^{+}=t_{j}$ and $t_{\ell}^{-}=t_{j+1}$. The
$t^{\pm}_{\ell}$ are such that the integration limits of
$u^{(\ell)}(t_{n})$ fit into the approximation interval: $t_{n} -
t^{\pm}_{\ell} \in I_{\ell}$. 
In Figure~\ref{fig:mosaicYM} we show how the solutions to the linear
ODEs $\bar y(t^{-}_{\ell},t^{+}_{\ell},\lambda^{\ell}_k)$ fill the
mosaic. In this example we have $B=3$ and the $t_j$ are given by the
non-equidistant sequence of time points indicated both along the
horizontal axis and the diagonal of the triangle. We can see here
that, for instance, when times greater than $t=3.2$ are reached, all
the ``past'' from $0$ to $2.11$ is stored in the solution values
$y(2.1,0,\lambda_k^{(\ell)})$, represented in
Figure~\ref{fig:mosaicYM} by the tall dark rectangle with basis
$[3.2, 4.1]$ for $\ell=3$, and in the adjacent incomplete white
rectangle to the right for $\ell=4$. In our example the values
$y(2.1,0,\lambda_k^{(3)})$ are used $8$ times to evaluate $u(t)$ for
$t\in [3.2,4.1]$. The filling of the mosaic is done bottom up in the
algorithm, advancing all the differential equations in every time
step, so that the algorithm can forget all the past values of the
function $g$, with the exception of those at $t_\ell^\pm$, which are
needed for the direct steps described below. In addition to the
current solution values $y(t,t_\ell^-,\lambda_k^{(\ell)})$ of the
differential equations at time $t=t_n$, also their values at a
splitting point $t_\ell^+$ need to be stored  until $t_\ell^+$ is
increased at a later step. Actually the algorithms stores three copies
of $y(t_{\ell}^{+},t_\ell^-,\lambda_k^{(\ell)})$.
 Pseudocodes for the organization of the
decomposition are given in the appendix.
\subsection{Direct steps}\label{Sec:dsteps}
The gaps $[t^{+}_{\ell+m},t^{-}_{\ell} ]$  between
 the enclosed blocks in Figure~\ref{fig:mosaicYM}
 are bridged using the values $\tilde{u}^{(\ell)}(t_n)$ whose
computation we describe next.
These direct steps compute
\begin{equation}
  \label{eq:dstep}
\tilde{u}^{(\ell)}(t_n) =\int_{t^+_{\ell+m}}^{t^-_{\ell}}
f(t_n-\tau)g(\tau)\, d\tau = \int_{t_j}^{t_{j+1}}
f(t_n-\tau)g(\tau)\, d\tau
\end{equation}
for some $j\in\{0,1,\dots,N-1\}$. On the
interval $[t_{j},t_{j+1}]$ we approximate $g(t)$ by a linear function:
$$
g(t)\approx \bar g(t)= g_j + \frac{g_{j+1}-g_j}{h_{j+1}}(t-t_j), \qquad
h_{j+1}= t_{j+1}-t_j,
$$
with $g_j=g(t_j)$, $j=0,1,\dots,N$.
(Here again, the approach would extend to polynomials of higher degree.)
Extending $\bar g(t)$ to $[0,t_{n}]$
we split \eqref{eq:dstep} in two terms
\begin{eqnarray*}
\lefteqn{ \int_{t_j}^{t_{j+1}} f(t_n-\tau)\,\bar g(\tau)\, d\tau =
    \int_{t_j}^{t_n}
    f(t_n-\tau)\,\bar g(\tau)\, d\tau - \int_{t_{j+1}}^{t_n}
    f(t_n-\tau)\,\bar g(\tau)\,
    d\tau} \hspace{1em}
  \\[5pt]
  &=&  \mathcal{L}^{-1}[F\cdot
  \mathcal{L}g(\cdot + t_j)](t_n-t_j) - \mathcal{L}^{-1}[F\cdot
   \mathcal{L}g(\cdot + t_{j+1})](t_n-t_{j+1})
  \\[5pt]
 &=&  \mathcal{L}^{-1}\bigg[ F_1 g_j +
  F_2\, \frac{g_{j+1}-g_j}{ h_{j+1}}
  \bigg](t_n-t_j)
  -\mathcal{L}^{-1}\bigg[F_1 \, g_{j+1} +
  F_2\, \frac{g_{j+1}-g_j}{ h_{j+1}} \bigg](t_n-t_{j+1}),
\end{eqnarray*}
where $F_1(s)=F(s)/s$ and $F_2(s)=F(s)/s^2$.
We approximate the
inverse Laplace transforms
$$
f_{1}(t) =
\bigl(\mathcal{L}^{-1}F_1\bigr)(t),
\qquad f_{2}(t) =  \bigl(\mathcal{L}^{-1}F_2\bigr)(t) \\
$$
at $t=t_n-t_j$ and $t=t_n-t_{j+1}$
using the numerical integration of the Laplace inversion formula
along the integration contours corresponding to the approximation
intervals $I_{\ell_{1}}$ and $I_{\ell_{2}}$ such that $t_n-t_{j+1}
\in I_{\ell_2}$ and $t_n-t_{j} \in I_{\ell_{1}}$.
The result of the direct step is then calculated forming linear
combinations
\begin{gather}
  \label{eq:linapprox}
  \begin{aligned}
    \int_{t_j}^{t_{j+1}} f(t_n-\tau)\,\bar{g}(\tau)\, d\tau  &=
    f_{1}(t_{n}-t_{j})\, g_{j} + f_{2}(t_{n}-t_{j}) \,\frac{g_{j+1} -
      g_{j}}{h_{j+1}}
    \\
    &  - f_{1}(t_{n}-t_{j+1})\, g_{j+1} - f_{2}(t_{n}-t_{j+1})\,
    \frac{g_{j+1} - g_{j}}{h_{j+1}}~.
  \end{aligned}
\end{gather}
This is also used for the terms closest to the diagonal,
$t_{j+1} = t_{n}$, where we note in addition that $f_1(0)=f_2(0)=0$.
\subsection{Complexity}
Given an arbitrary sequence of time points $0<t_1<\ldots<t_N=T$ with
the minimum step size $h_* = \min_j (t_{j+1}-t_j)$, the above
algorithm computes
\begin{equation}\label{eq:convollin}
\int_0^t f(t-\tau)\, \bar g(\tau)\, d\tau, \quad\hbox{ for }\
t=t_1,\dots,t_N,
\end{equation}
(with $\bar g$ the piecewise linear interpolant of $g$) up to an
error $\veps$ using
\begin{itemize}
\item $L=O(\log\frac T{h_*})$ hyperbolas with
\item $2K+1=O(\log\frac1\veps)$ quadrature points on each hyperbola.
\end{itemize}
The algorithm thus requires $(2K+1)L$ evaluations of the Laplace
transform $F(s)$ at the quadrature points and solves $(2K+1)L$
differential equations $y'=\lambda_k^{(\ell)}y + \bar g$. As the
algorithm proceeds, only three solution values need to be stored for
each of these differential equations. In addition, at most $2L$ values
of $g$ need to be
kept in memory for the direct steps. In total, the active memory
requirements are $O(LK)=O(\log\frac T{h_*}\log\frac1\veps)$ vectors
of the dimension of $g$. The total operation count is
$O(NLK)=O(N\log\frac T{h_*}\log\frac1\veps)$. For the
variable-step-size algorithm we thus obtain the complexity
characteristics as stated in the introduction for the fast and
oblivious fixed-step-size algorithm, with $\log N$ now replaced by
$\log\frac T{h_*}$.

\subsection{Adaptivity based on controlling the interpolation error}
\label{Sec:strategy_interpolationerror} There are two sources of
error in the algorithm. The first one is the discretization of the
contour integral, which is well controlled. The second one is the
piecewise linear interpolation of $g$ by $\bar g$. Ignoring the
error from discretizing the contour integrals, the algorithm thus
computes (\ref{eq:convollin}) instead of (\ref{eq:convolution}). We
control the error in $g$, which is bounded by
$$
\| \bar g(t) - g(t) \| \le
{\textstyle\frac18}\, h_n^2\,\max_{t_{n-1}\le\tau\le t_n} \| g''(\tau) \|
 \quad\ \hbox{ for}\quad t_{n-1}\le t \le t_{n}.
$$
Given a tolerance {\it Tol}, we propose the new step-size
$h_{n+1}$ according to the criterion
\begin{equation}\label{eq:new-step}
C h_{n+1}^2 \gamma_n'' = 0.8\cdot \mbox{\it Tol}\,,
\end{equation}
where $\gamma_n''=\| {\widetilde g}''(t_n) \|$ with the quadratic interpolant
$\widetilde g$ to $g$ at $t_{n-2}, t_{n-1}, t_n$. The constant $C$ is
chosen as
$
C \approx \frac18 \int_0^T |f(t)|\, dt.
$
Additionally the step-size is restricted to fulfill
$1/2 h_{n} < h_{n+1} < 2h_{n}$.
The proposed step size $h_{n+1}$ is tested by
\begin{equation}\label{eq:step-test}
C h_{n+1}^2 \gamma_{n+1}'' \le \mbox{\it Tol}\,,
\end{equation}
where the new value $g(t_n+h_{n+1})$ is used in the
computation of $\gamma_{n+1}''$.
If this condition is satisfied, then $h_{n+1}$ is accepted and we set
$t_{n+1}=t_n+h_{n+1}$,
else we repeat the test
with a reduced step $h_{n+1}$ determined from (\ref{eq:new-step}) with
$\gamma_{n+1}''$ in place of $\gamma_n''$. If necessary,
this procedure is repeated until (\ref{eq:step-test}) is satisfied.

In the following sections we give three examples to
show the performance of the algorithm, using also other strategies
for controlling
the step size in the time integration of integro-differential equations.
However, we point out that the above fast algorithm is
independent of the particular step size selection strategy.
The step size control is just the way we
generate the sequence of time points. The minimum step size need not be
specified {\it a priori}, and at time $t_n$ the future time points
$t_{n+1},t_{n+2},\dots$ need not yet be known in the algorithm.

For the three examples provided, we ran the algorithm with basis
$B=5$, which gives approximation intervals of the type
$[t_0,25t_0]$.

\section{A blow-up example originating from a Schr\"odinger equation with
concentrated nonlinearity}
Adami \& Teta \cite{AdaT99,AdaT01} consider a nonlinear Schr\"odinger
equation with nonlinearity concentrated at $x=0$,
\begin{equation}\label{eq:conc-schroed}
i\frac{\partial\psi}{\partial t} = -\Delta\psi +
\gamma |\psi|^{2\sigma}\psi\cdot \delta_{x=0}\,,
\qquad x\in\Rset,\ t>0,
\end{equation}
with initial data $\psi(x,0)=\psi_0(x)$ for $x\in\Rset$.
The equation can be given a rigorous formulation as an integral
equation. With Duhamel's principle the wave function
$\psi(x,t)$ can be expressed as the sum
of $\bigl(U(t)\psi_0\bigr)(x)$, that is,
the solution at position $x$ and time $t$
of the free Schr\"odinger equation with initial data $\psi_0$,
and of a convolution with $|z|^{2\sigma}z$, where $z(\tau):=\psi(0,\tau)$
for $0\le\tau\le t$. The function $z$ is the solution of a
nonlinear, complex Abel-type integral equation
\begin{equation}\label{eq:abel}
z(t) + \gamma\, \frac{\sqrt{i}}2 \int_0^t \frac1{\sqrt{\pi(t-\tau)}}\,
|z(\tau)|^{2\sigma} z(\tau)\, d\tau = a(t)\,,\qquad t\ge 0,
\end{equation}
where $a(t)=\bigl(U(t)\psi_0\bigr)(0)$ is the solution at $x=0$
of the free Schr\"odinger equation. We present the results of numerical
experiments in a situation where the solution is known to have
finite-time blow-up. We choose $\sigma=1$, and
$$
a(t)=\frac1{\pi^{1/4}}\frac1{\sqrt{1+i2t}}\,,
$$
which corresponds
to a Gaussian wave-packet $\psi_0(x)=\pi^{-1/4}e^{-x^2/2}$ as initial data.

\begin{figure}[!h]
 \centering
  \includegraphics[width=0.48\textwidth]{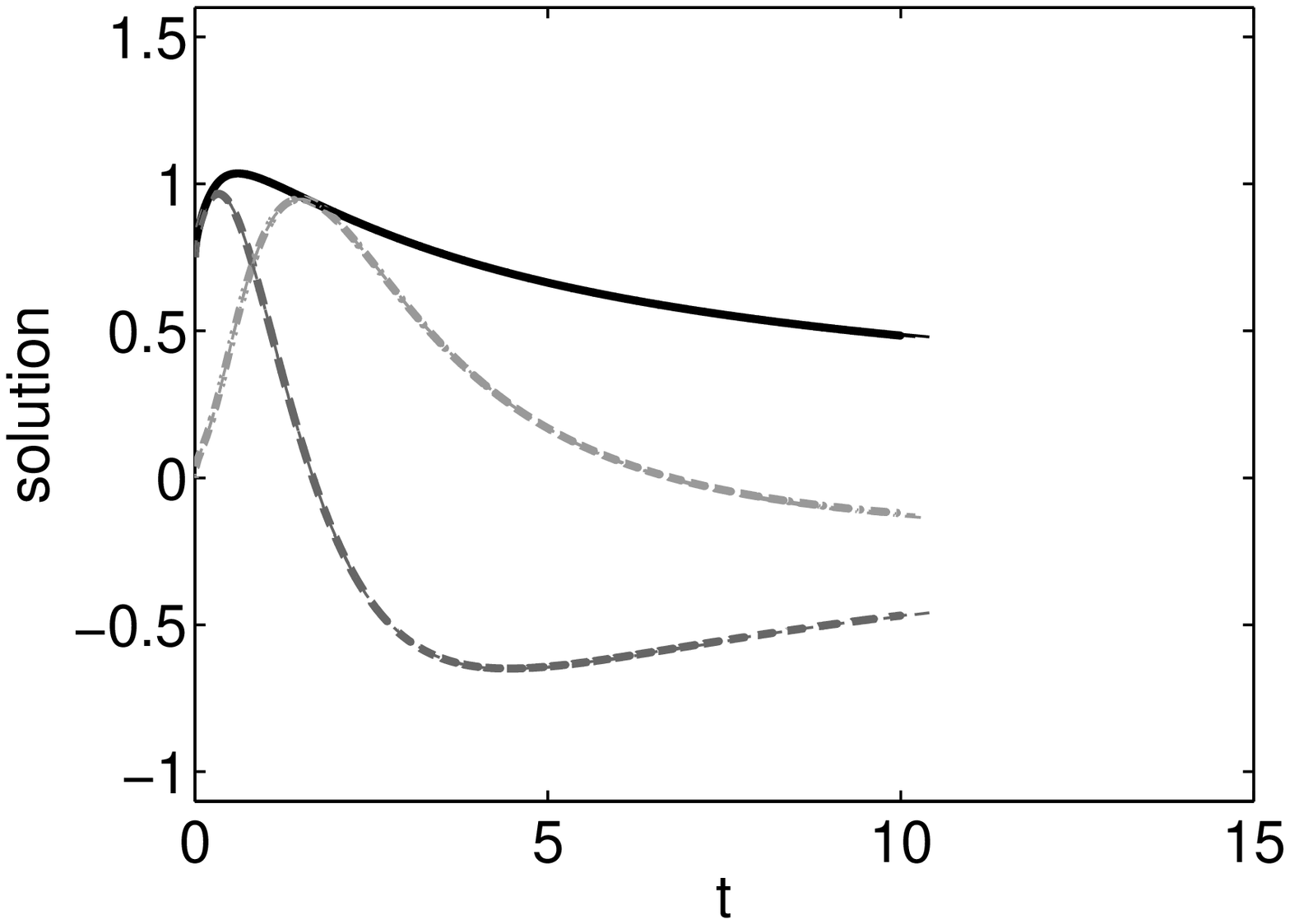}
  \includegraphics[width=0.48\textwidth]{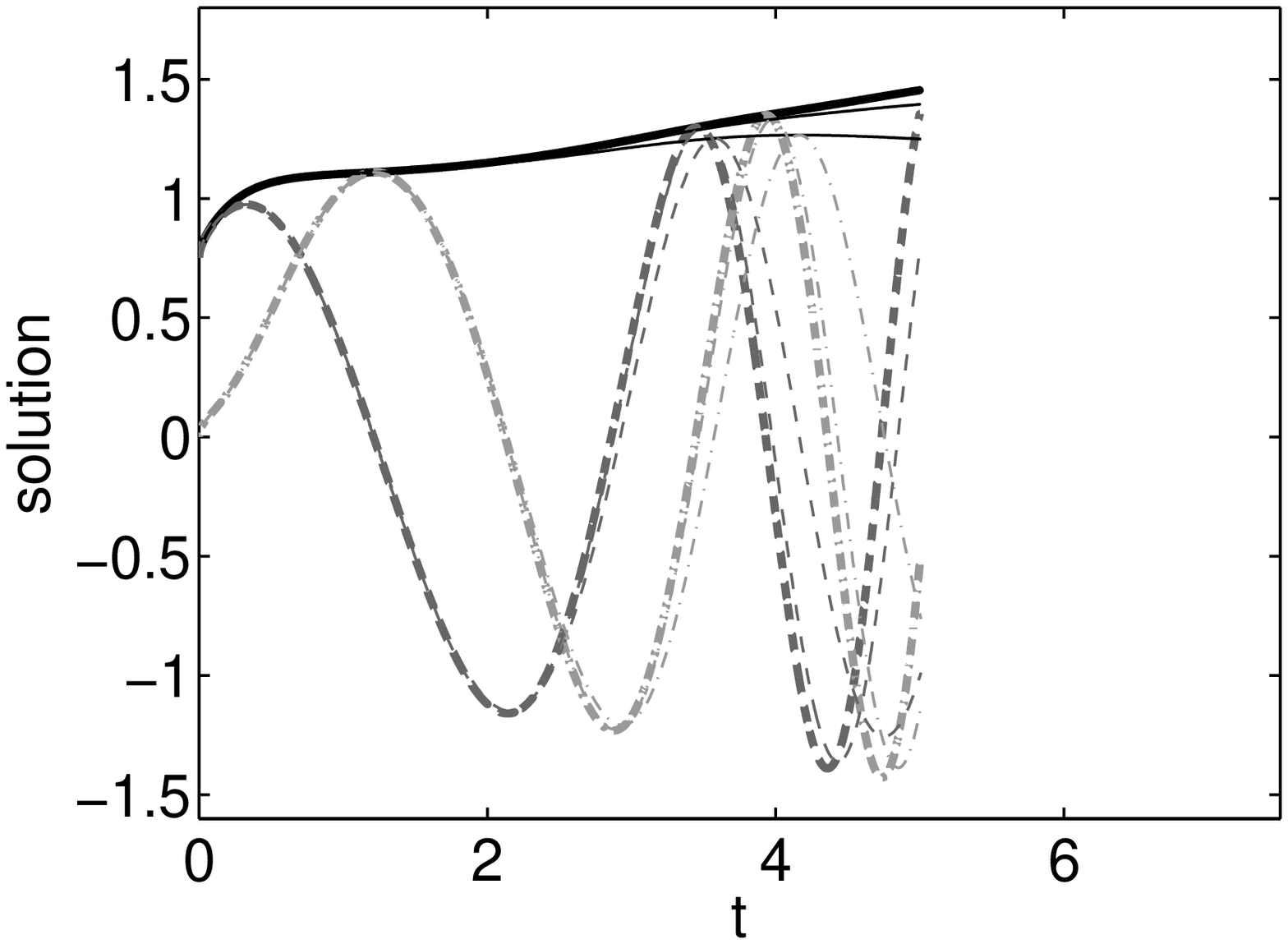} \\
  \includegraphics[width=0.48\textwidth]{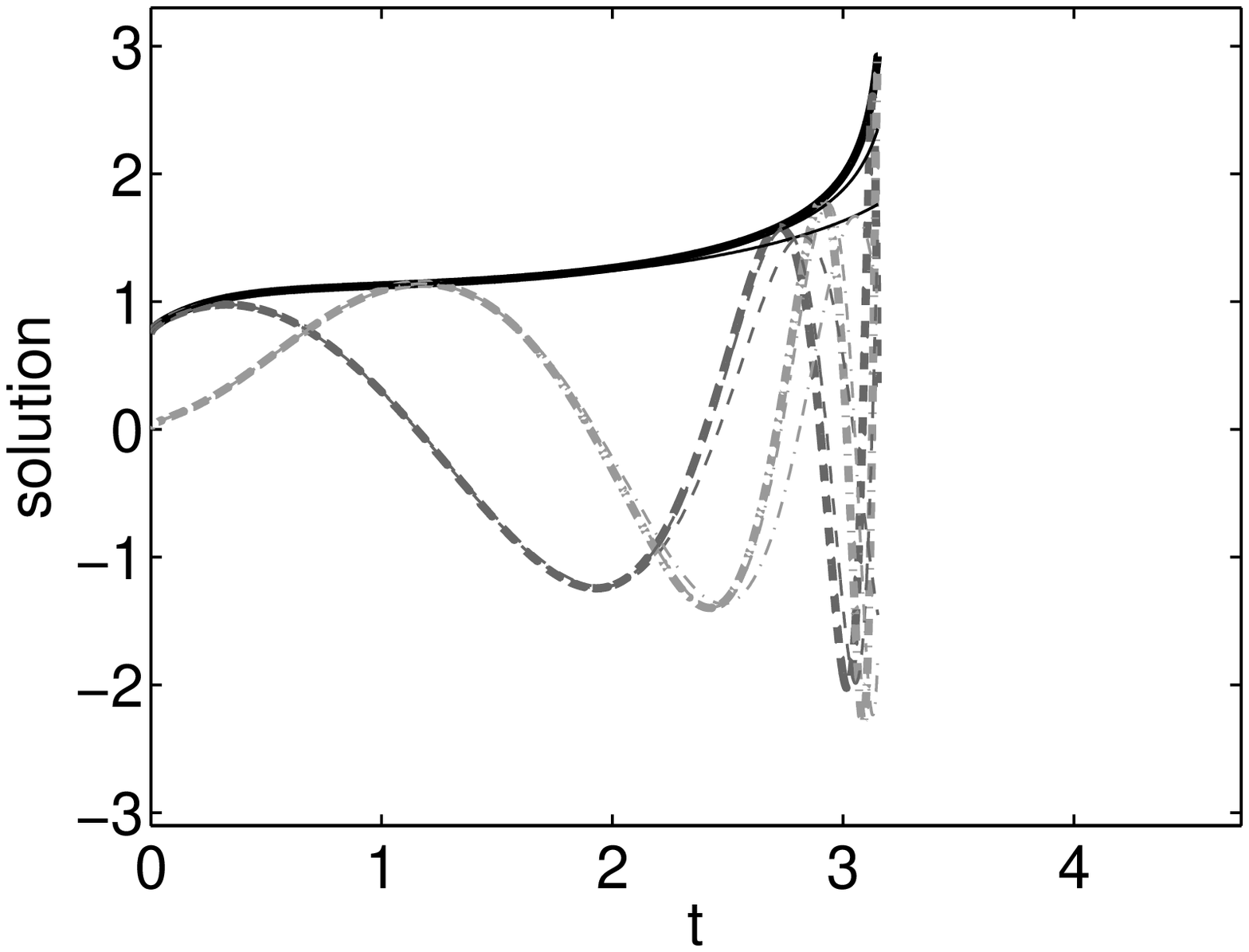}
  \includegraphics[width=0.48\textwidth]{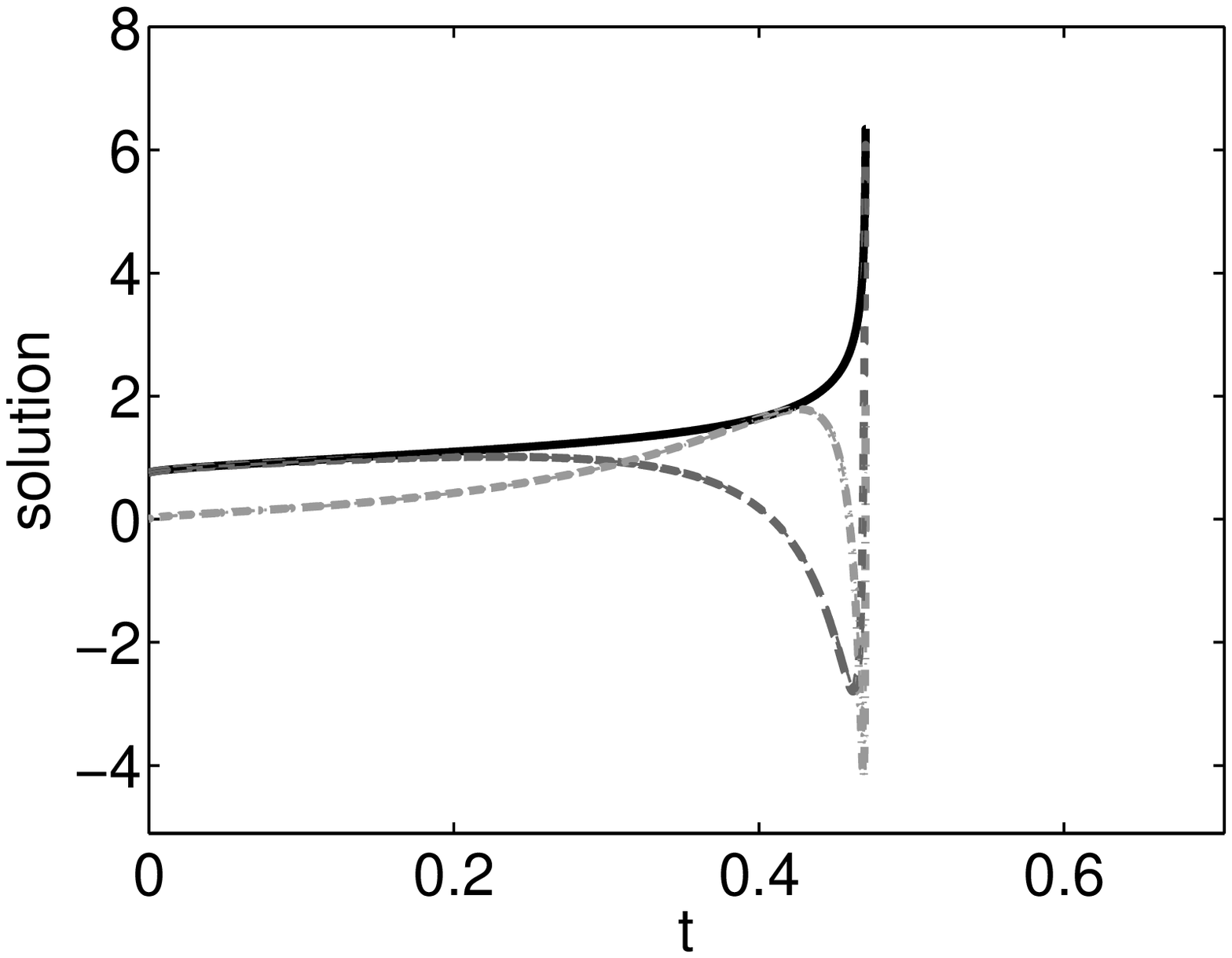}
  \caption{Real and imaginary part (light and dark gray) and modulus (black) of the solution $z$
    for $\gamma= -2, -2.05, -2.06 -2.5$ and for different tolerances.}
  \label{fig:sol_blowup}
\end{figure}
For four different values of $\gamma$, $\gamma = -2, -2.05, -2.06, -2.5$
the evolution of the solution is shown in Fig~\ref{fig:sol_blowup} for
different tolerances. The thick lines correspond to solutions obtained with
a tolerance of $10^{-7}$. The other two lines correspond to tolerances of
$10^{-3}$ and $2 \cdot 10^{-4}$. Whereas for $\gamma =-2,\, -2.5$
one cannot distinguish the different tolerances in the plot, for $\gamma =-2.05,\, -2.06$
it is clearly visible, that choosing a too low tolerance will produce a wrong result.

\begin{figure}[!h]
  \centering
  \includegraphics[width=0.45\textwidth]{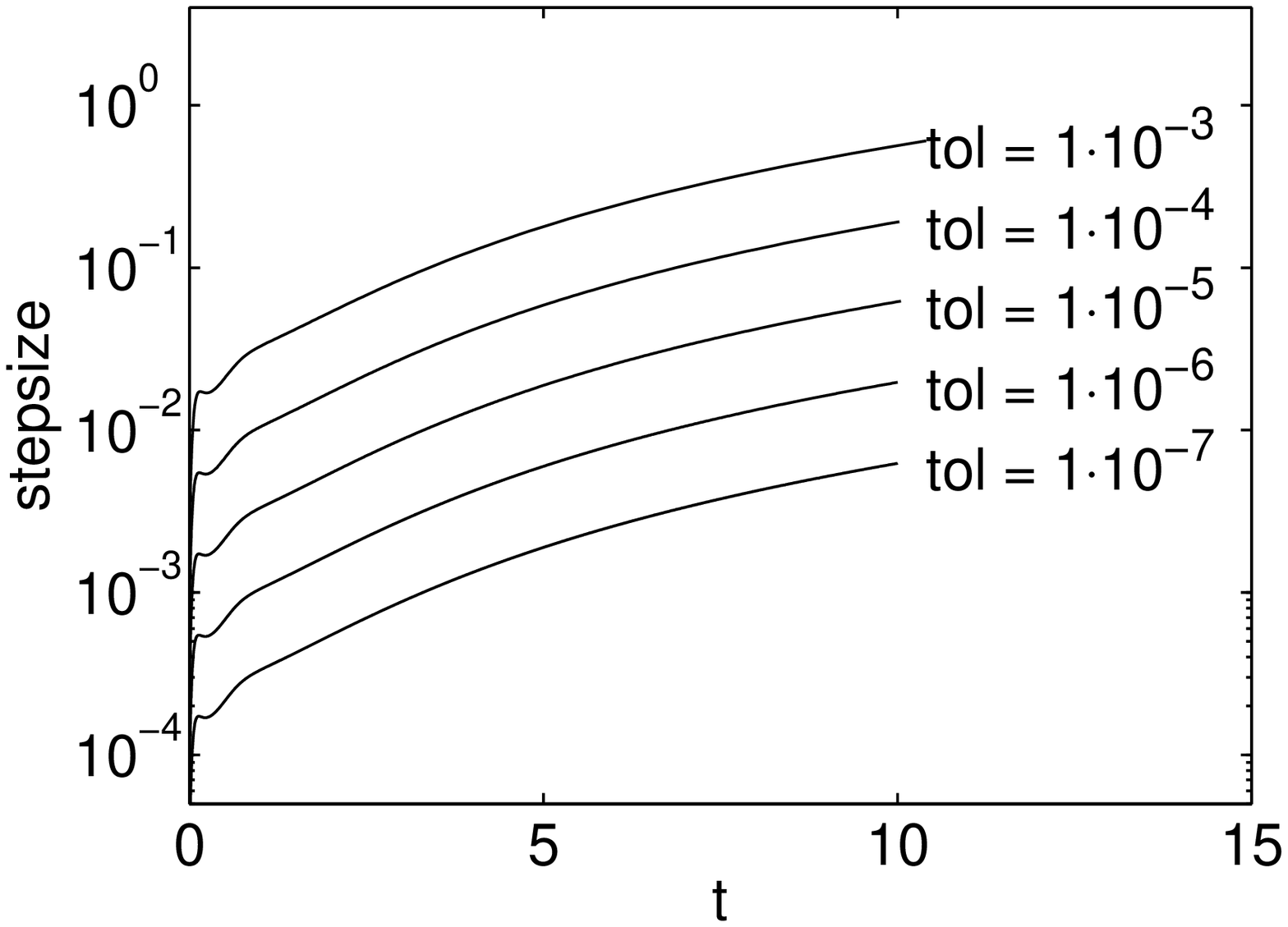}
  \includegraphics[width=0.45\textwidth]{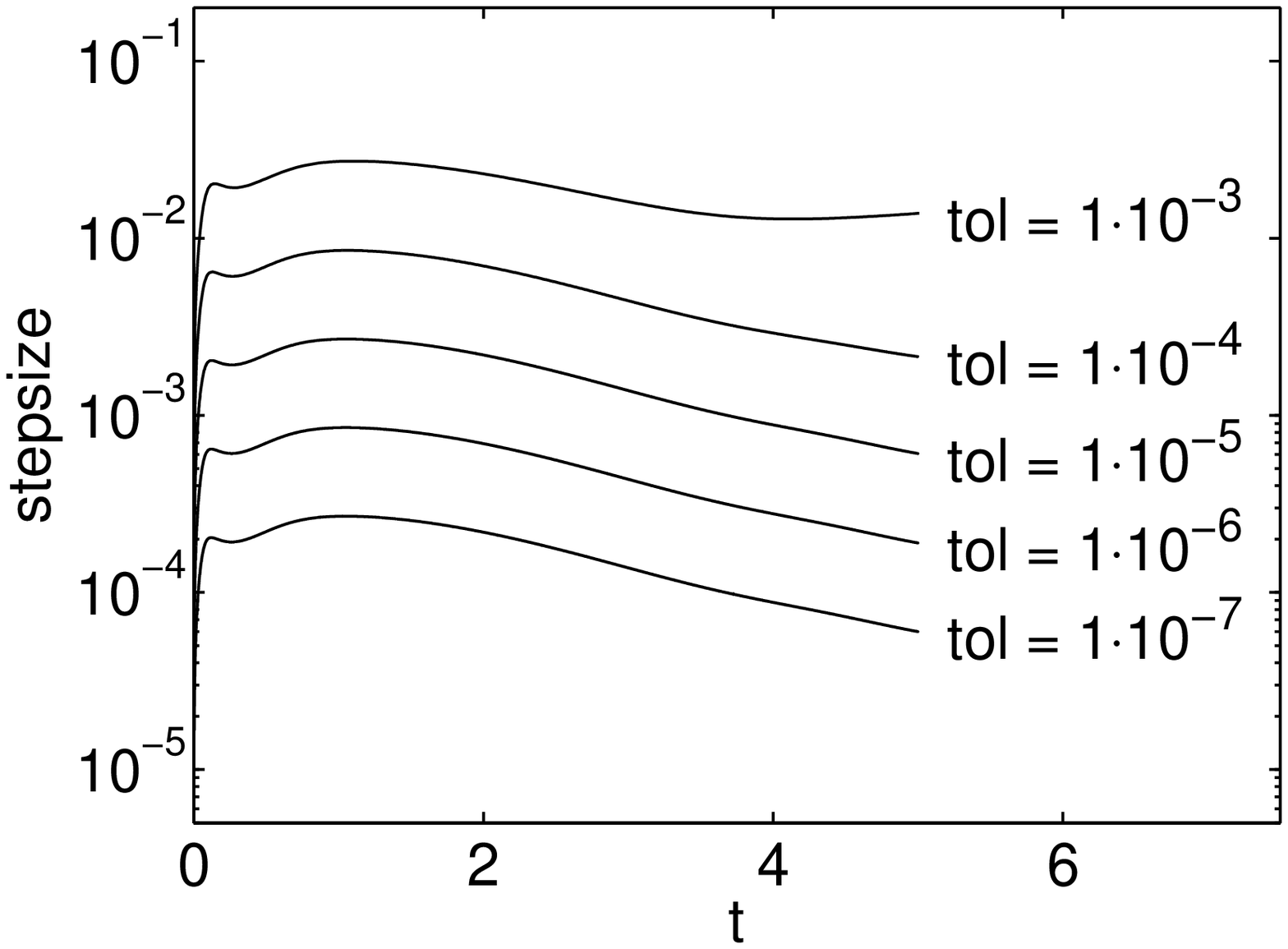}\\
  \includegraphics[width=0.45\textwidth]{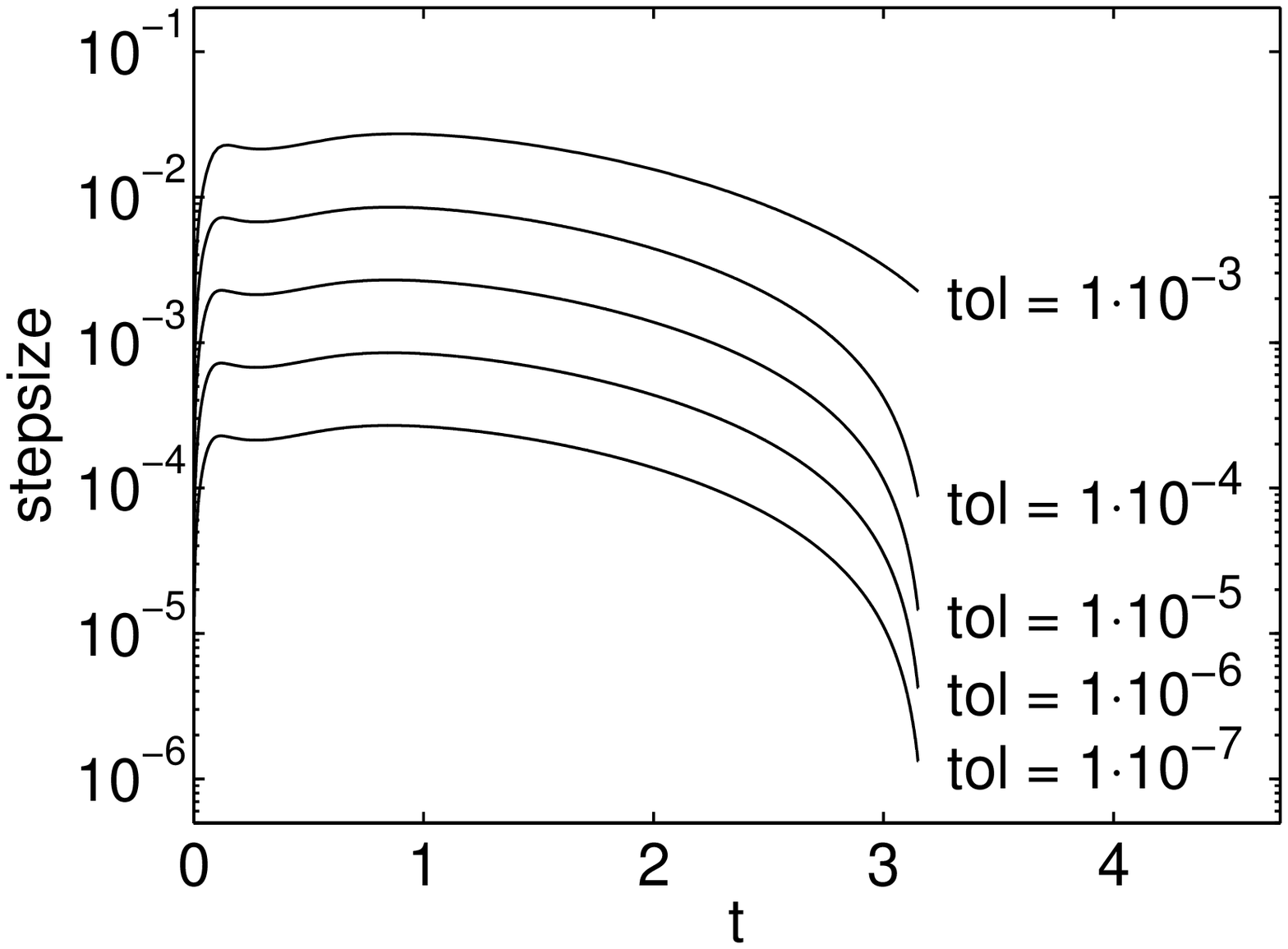}
  \includegraphics[width=0.45\textwidth]{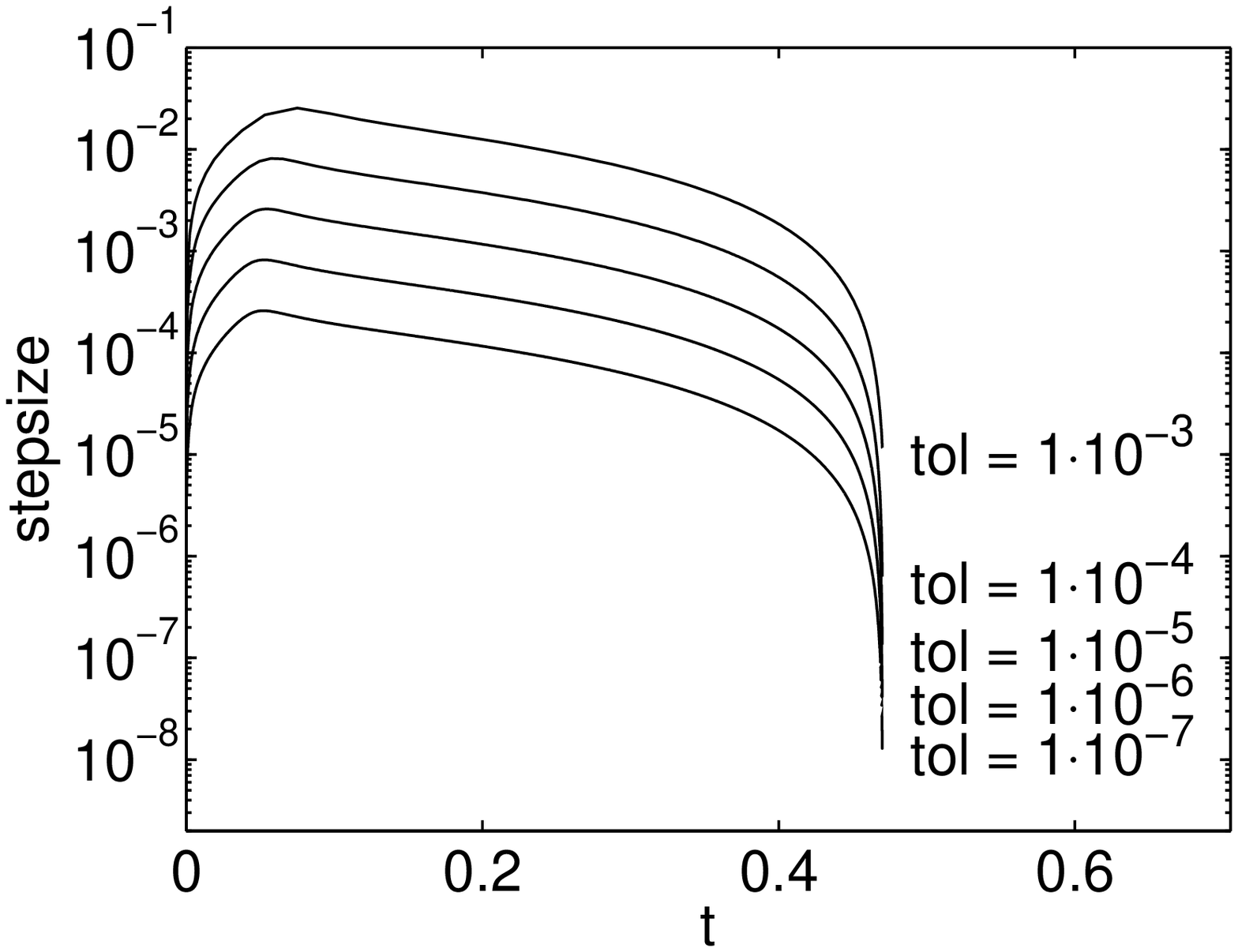}
  \caption{Evolution of the step-size for $\gamma= -2, -2.05, -2.06, -2.5$
    and different tolerances}
  \label{fig:stepsize_blowup}
\end{figure}
In Fig.~\ref{fig:stepsize_blowup} the step size is plotted over $t$. Clearly the
adaptivity pays off to resolve the blow-up. In case $\gamma=-2.5$ we observe for
a tolerance of $10^{-7}$ step-sizes ranging from $10^{-8}$ to $10^{-4}$.

\begin{figure}[!h]
  \centering
  \includegraphics[width=0.45\textwidth]{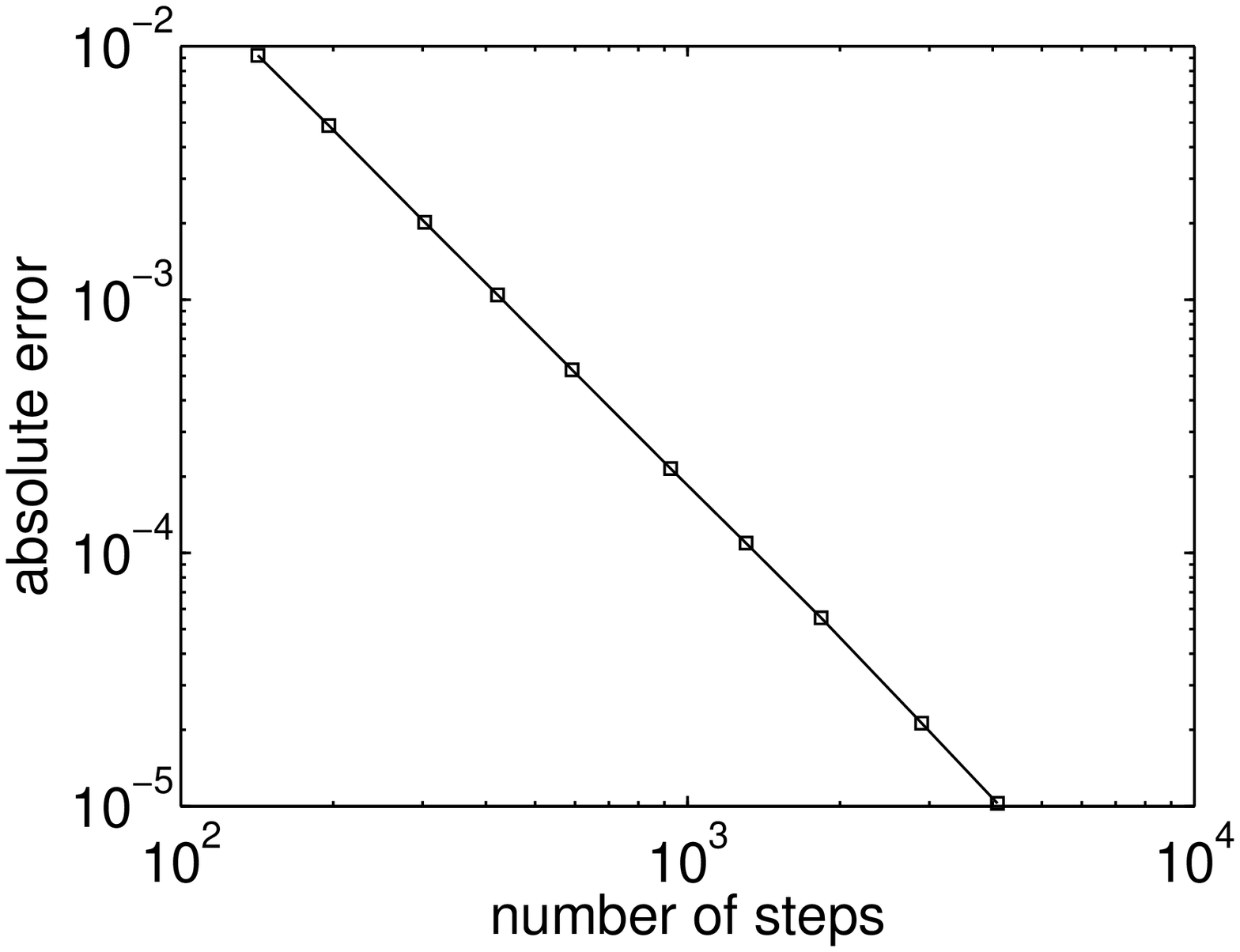}
  \includegraphics[width=0.45\textwidth]{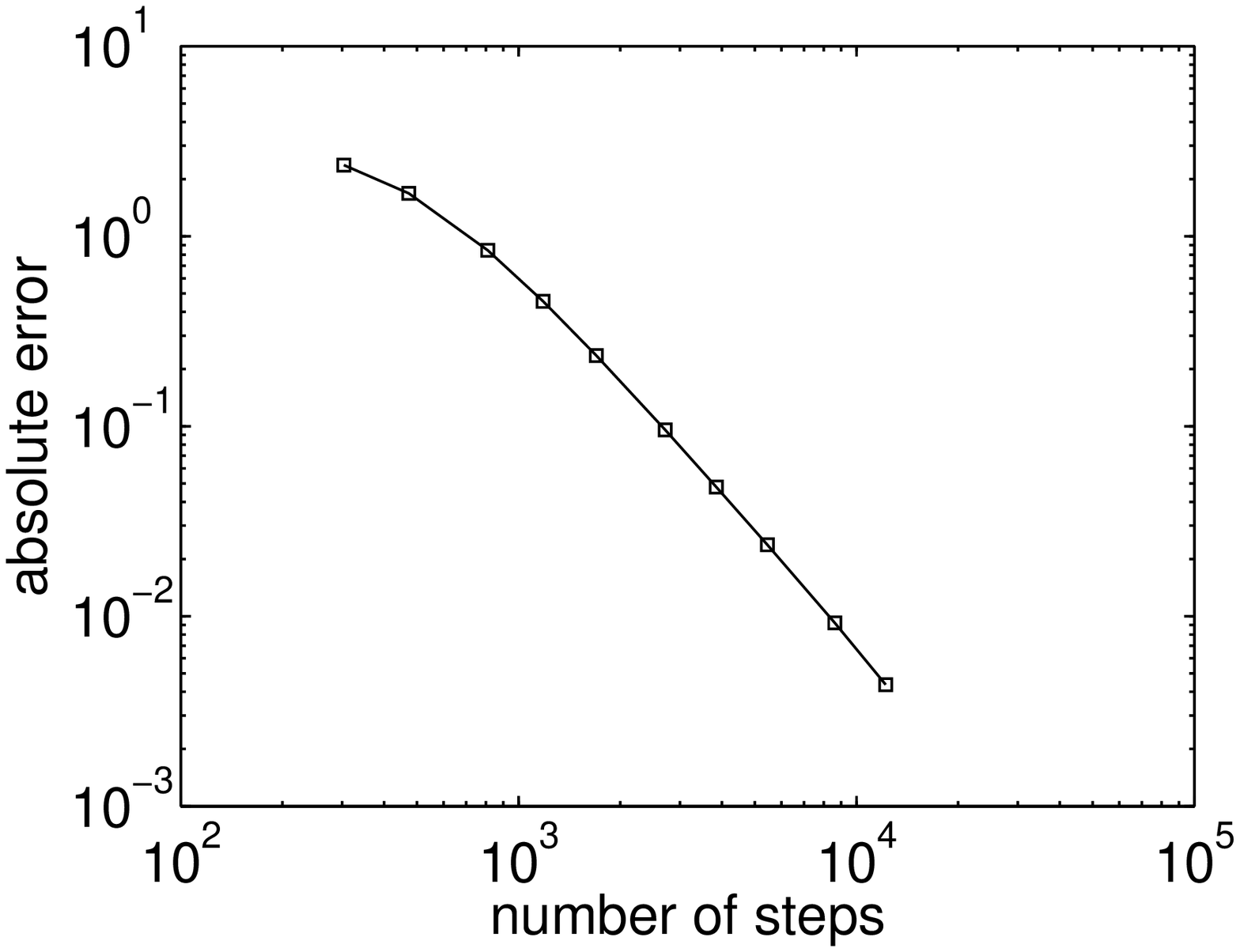}\\
  \includegraphics[width=0.45\textwidth]{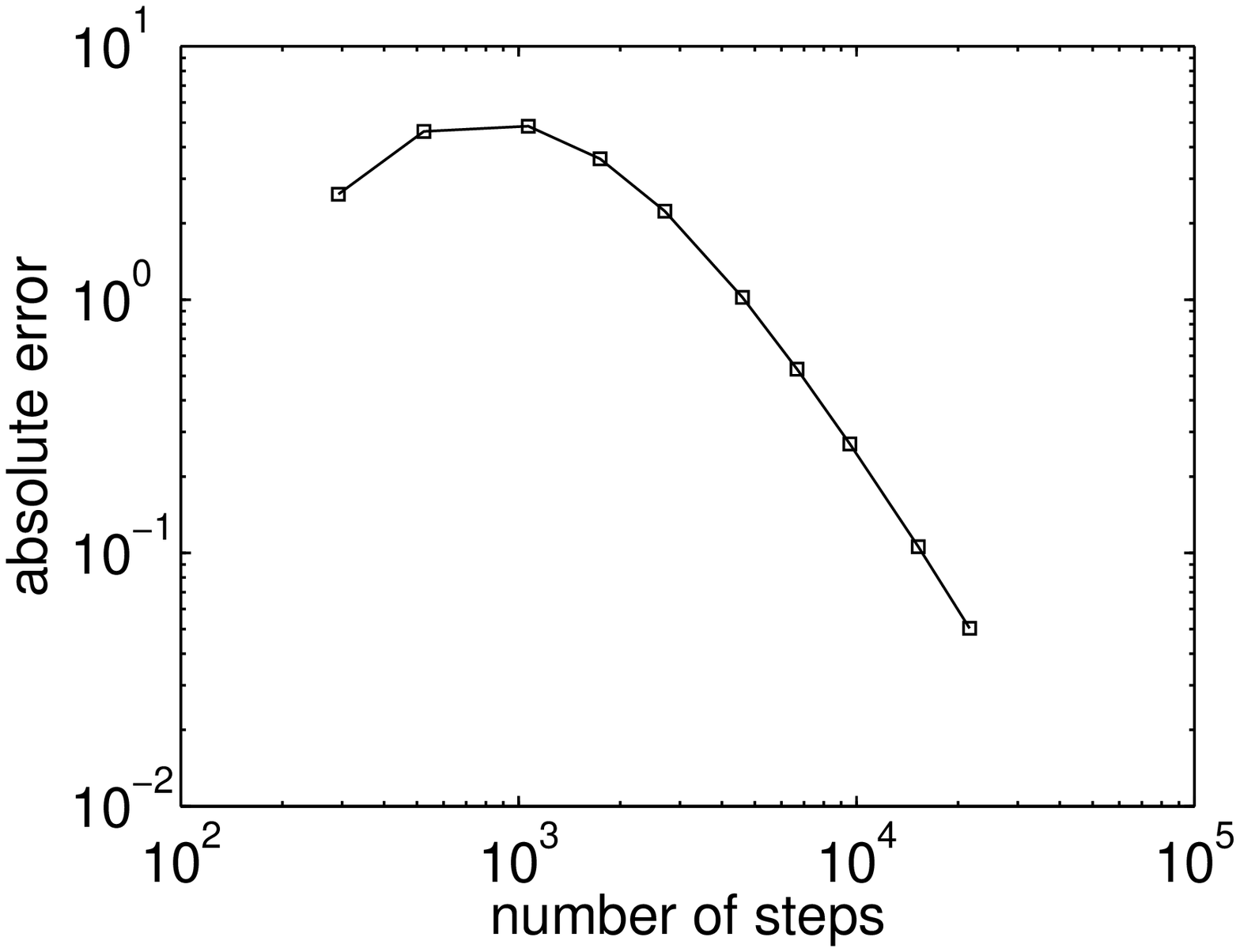}
  \includegraphics[width=0.45\textwidth]{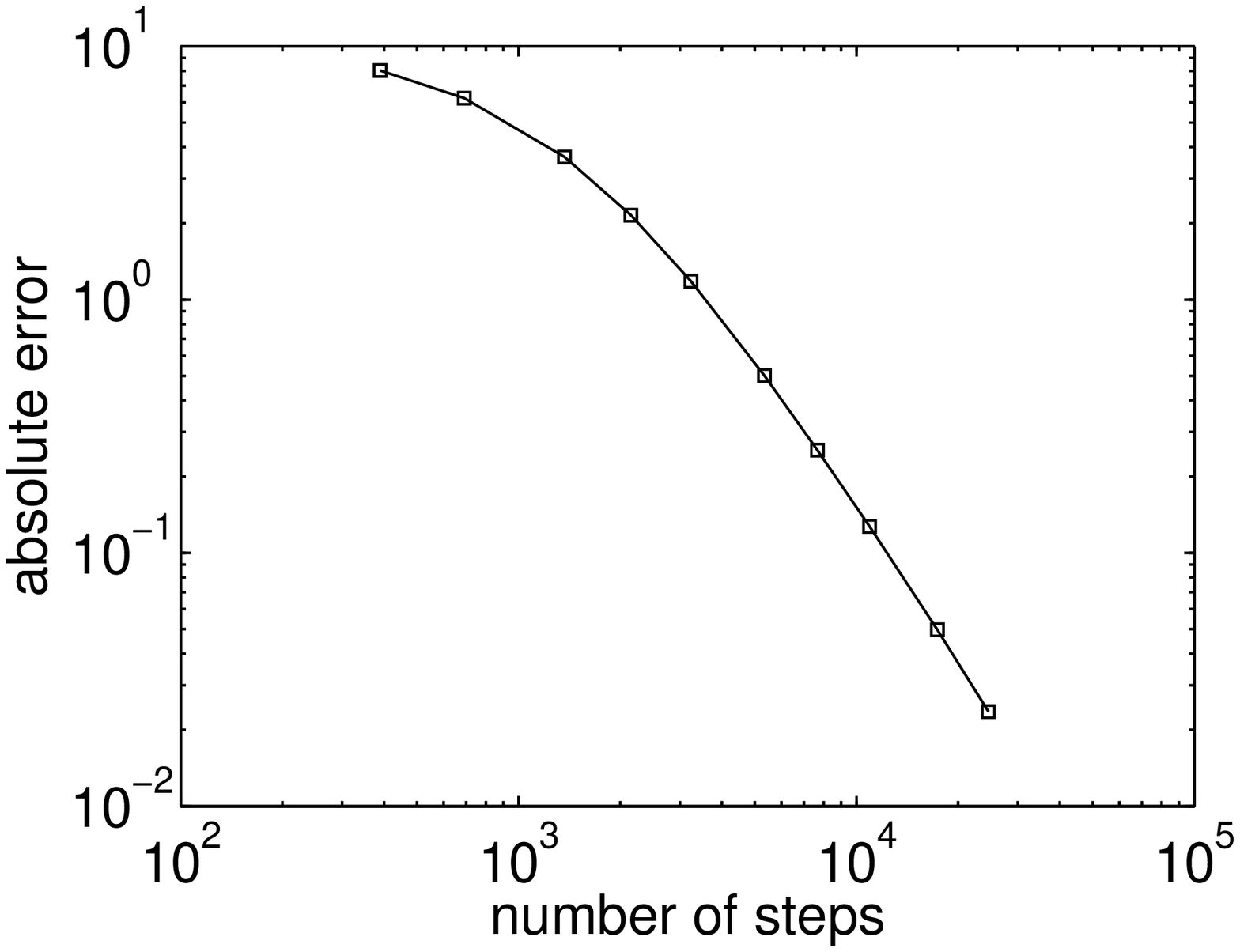}
  \caption{Error versus the number of steps,
    at $t=10$ for $\gamma=-2$, at $t=5$ for $\gamma=-2.05$, at $t=3.15$
    for $\gamma=-2.06$, and at $t=0.47$ for $\gamma=-2.5$}
  \label{fig:err_ns_blowup}
\end{figure}
Fig.~\ref{fig:err_ns_blowup} displays the absolute error at the
final time, obtained with tolerances $1\cdot 10^{-3}, 5\cdot
10^{-4}, 2\cdot 10^{-4}, 1\cdot 10^{-4}, 5\cdot 10^{-5}, 2\cdot
10^{-5}, 1\cdot 10^{-5}, 5\cdot 10^{-6}, 2\cdot 10^{-6}, 1\cdot
10^{-6}$ versus the total number of steps. The error is measured
against a reference solution obtained with tolerance $10^{-7}$. The
numerical inversion of the Laplace transform is performed as
explained in Section~\ref{Sec:InvLaplace}, with $a=0.8$, $d=0.7$,
and $K=50$, which gives $C_1=6.567$ and $C_2=0.066$ (see
Theorem~\ref{thm:err}). For larger tolerances, good results can be
obtained with a smaller $K$, say $K=25$. Taking $K=40$, we
only get small oscillations in the stepsize for the smallest
tolerance, $\mathtt{tol}=10^{-7}$, in
Figure~\ref{fig:stepsize_blowup}, and no visible changes for less
stringent tolerances.

\section{Chemical reaction kinetics with inhibited diffusion}
We consider three molecular species $A$, $B$ and $C$, reacting as
\begin{eqnarray}\label{eq:rd_ref}
A + B &\stackrel{k_{1}}{\longrightarrow}& C
\nonumber\\
C &\stackrel{k_{2}}{\longrightarrow}& A+B
\nonumber\\
C &\stackrel{k_{3}}{\longrightarrow}& A+P~,
\end{eqnarray}
$P$ being the resulting product.
The diffusion of each of the species is anomalous.
So we obtain a reaction diffusion equation with a memory term.
A model like this, with three species and
three reactions was considered in~\cite{CaoBA06}.
However we have chosen to follow~\cite{YuLi04} and to associate
a memory with the reaction term. Thus to model this process the following
system of integro-differential equations is considered:
\begin{gather}
\label{cauchy_rd}
\begin{aligned}
  \dot{u}_{1} =& \partial_t \partial_t^{-\alpha} \big( K \Delta u_1 - k_{1} u_1u_2 +
  (k_{2}+k_{3})u_3 \big)
  \\
  \dot{u}_{2} =& \partial_t \partial_t^{-\alpha} \big( K \Delta u_2 - k_{1} u_1u_2 +
  k_{2}u_3 \big)
  \\
  \dot{u}_{3} =& \partial_t \partial_t^{-\alpha} \big( K \Delta u_3 + k_{1} u_1u_2 -
  (k_{2}+k_{3})u_3\big),
\end{aligned}
\end{gather}
where $\Delta = \partial_{xx}$ is the 1D Laplacian with periodic
boundary condition on $[-5,5]$ and $\partial_t^{-\alpha}$ denotes
the fractional integral of order $0 < \alpha < 1$, given by the
Riemann Liouville operator
\begin{equation}\label{fracint}
 \partial_t^{-\alpha} g(t) =
\frac{1}{\Gamma(\alpha)} \int_{0}^{t}(t-\tau)^{\alpha-1} g(\tau) \,
d\tau,
\end{equation}
for $0< \alpha <1$.
Integrating in time system~\eqref{cauchy_rd}, we get the
integro-differential equation
\begin{gather}
  \label{volterraint_rd}
  \begin{aligned}
  u_{1}(t)- u_1(0) =& \partial_t^{-\alpha} \big[ K\Delta u_1(t) - k_{1} u_1(t) u_2(t) + (k_{2}+k_{3})u_3(t)
  \big]  \\
  u_{2}(t)- u_2(0) =& \partial_t^{-\alpha} \big[ K\Delta u_2(t) - k_{1} u_1(t) u_2(t) + k_{2}u_3(t)
  \big]
 \\
  u_{3}(t)- u_3(0) =& \partial_t^{-\alpha} \big[ K\Delta u_3(t) + k_{1} u_1(t) u_2(t) - (k_{2}+k_{3})u_3(t) \big].
  \end{aligned}
\end{gather}

In this situation we have in the convolution terms the weakly
singular kernel $f(t)=t^{\alpha-1}/\Gamma(\alpha)$, with Laplace
transform $F(s)=s^{-\alpha}$. We approximate the solution to
\eqref{volterraint_rd} by using an adaptive strategy similar to the
one explained in Section~\ref{Sec:strategy_interpolationerror} but
replacing criterion~\eqref{eq:new-step} by
$$
  C h_{n+1}^2 \gamma'_n = 0.8 \cdot Tol,
$$
and the test~\eqref{eq:step-test} by
$$
C h_{n+1}^2 \gamma_{n+1}' \le Tol\,,
$$
where $\gamma_n'=\| \widetilde{g}'(t_n) \|$, with
$\widetilde{g}$ the linear interpolant of $g$ at $t_{n-1}$ and
$t_n$. Our choice for the different parameters is $K = 0.5$, $k_1 =
1$, $k_2=2$ and $k_3=3$, and we integrate up to $T=30$. We fix
$\alpha=0.5$ and consider smoothed step-like functions as the
initial data.

Setting ${\bf u}=[u_1,u_2,u_3]^T$, $I_3$ the $3\times 3$ identity
matrix, and following the notation introduced in
Section~\ref{Sec:dsteps} for the direct steps of the algorithm, the
discrete equation approximating \eqref{volterraint_rd} is
\begin{gather}
  \label{eq:disu_rd_int}
  \begin{aligned}
    & \Big( I_3 \otimes I - \frac{f_2(h_n)}{h_n}\, ( K \, I_3 \otimes
    S + R)\Big){\bf u}^n
    \\
    &  = \Big( f_1(h_n)- \frac{f_2(h_n)}{h_n} \Big)\big( K
    \, I_3 \otimes S + R \big) {\bf u}^{n-1} +
    \,\, k_1 f_1(h_n)\, {\bf e} \otimes u_1^{n-1} u_2^{n-1} \\
    & + \,\, \frac 1{\Gamma(\alpha)} \int_0^{t_{n-1}}
    \frac {g(\bar{{\bf u}}(\tau))}{(t-\tau)^{1-\alpha}}\,d\tau +
    \mathbf{u}_0,
  \end{aligned}
\end{gather}
where $\bar{\bf u}$ denotes the piecewise linear interpolant of
${\bf u}$ at times $t_0,t_1,\dots,t_n$ and, for $M$ nodes in the
spatial discretization and ${\bf v}$ a column vector of length $3M$,
we define
$$
g({\bf v}) = \big( K \, I_3 \otimes S + R \big) {\bf v} + k_1 \big(
{\bf e} \otimes v_1 v_2 \big), \qquad {\bf e}=[-1, -1, 1]^T,
$$
with $S$ the second order finite difference approximation to
$\partial_{xx}$ with periodic boundary conditions and $R$ the $3M
\times 3M$ matrix
$$
R = \left(
\begin{array}{ccc}
0 & 0 & (k_{2}+k_{3})I_{M}
\\
0 & 0 &  k_{2}I_{M}
\\
0 & 0 & -(k_{2}+k_{3})I_{M}
\end{array}
\right),
$$
where $I_M$ is the $M\times M$ identity matrix.

\begin{figure}
  \centering
  \includegraphics[width=0.43\textwidth, height=0.2\textheight]{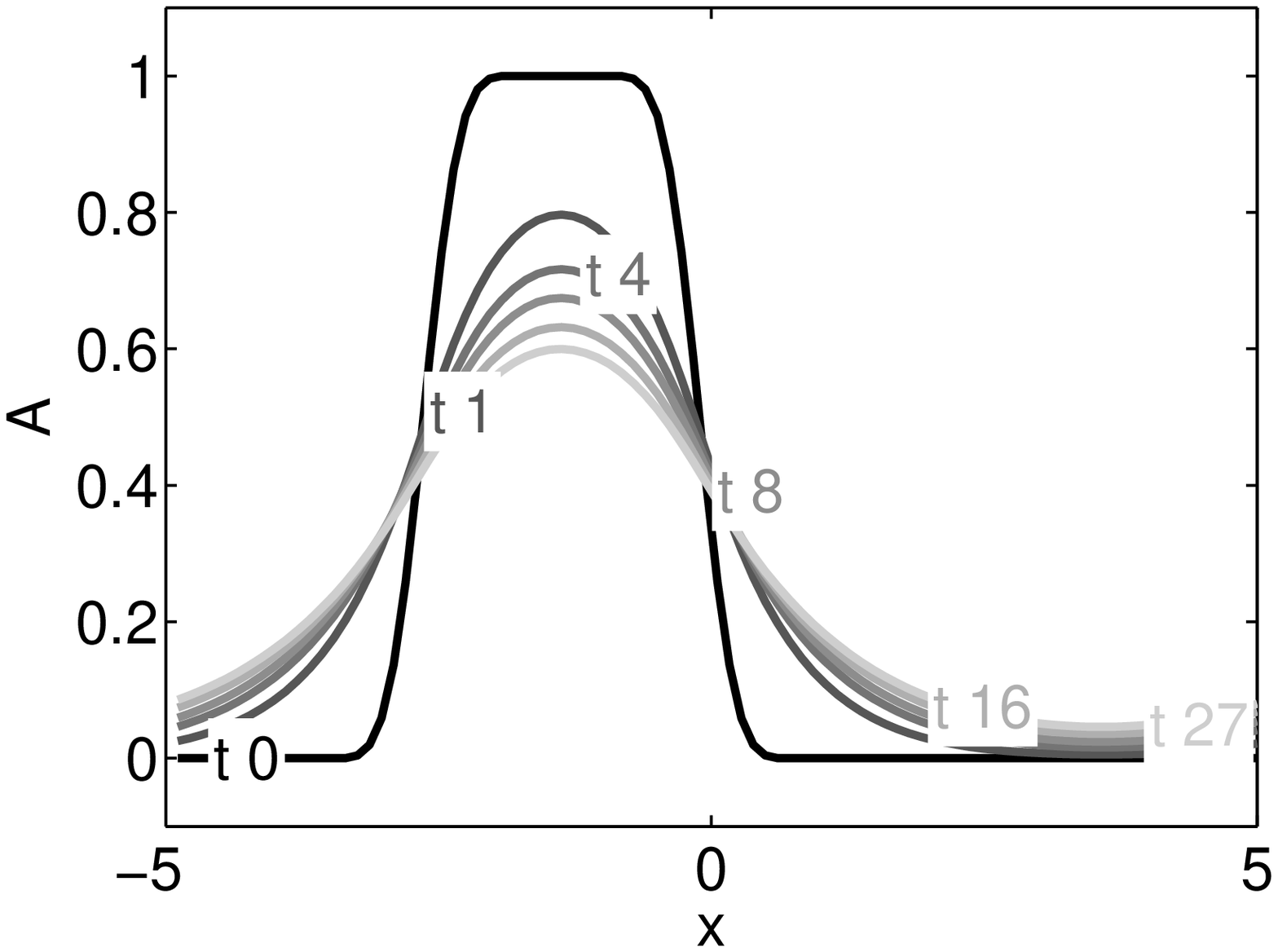}
  \includegraphics[width=0.43\textwidth, height=0.2\textheight]{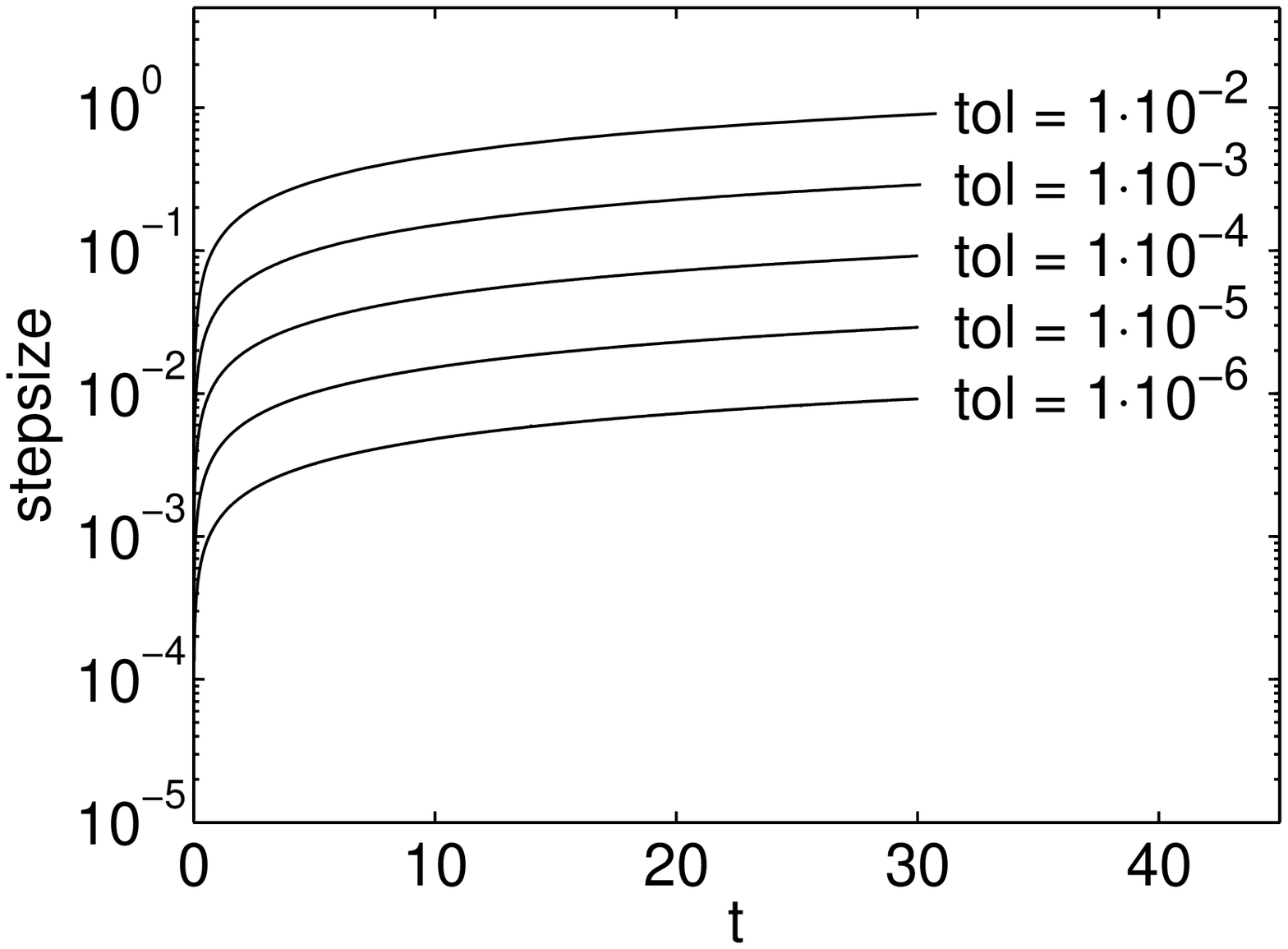}
  \\
  \includegraphics[width=0.43\textwidth, height=0.2\textheight]{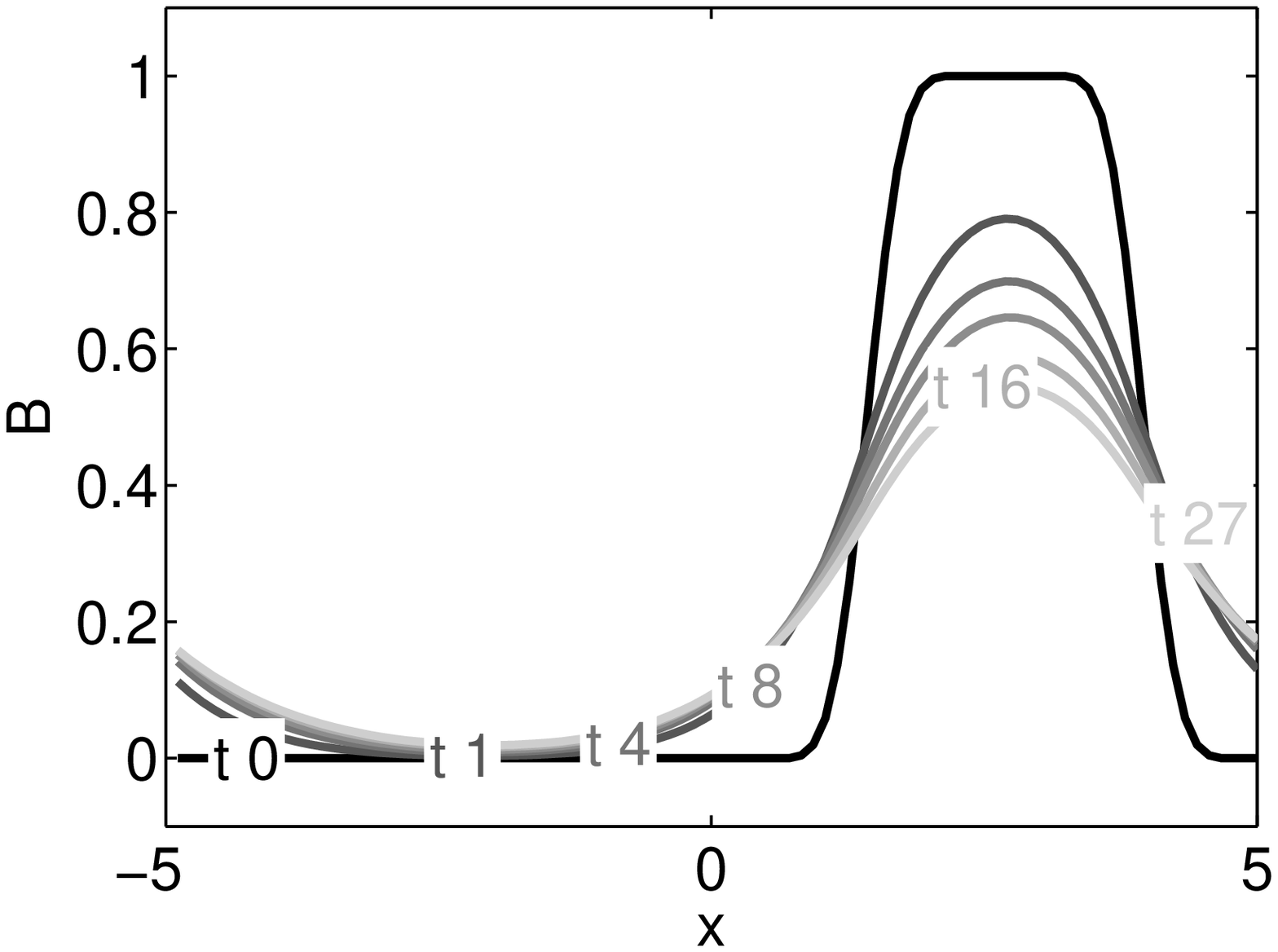}
  \includegraphics[width=0.43\textwidth, height=0.21\textheight]{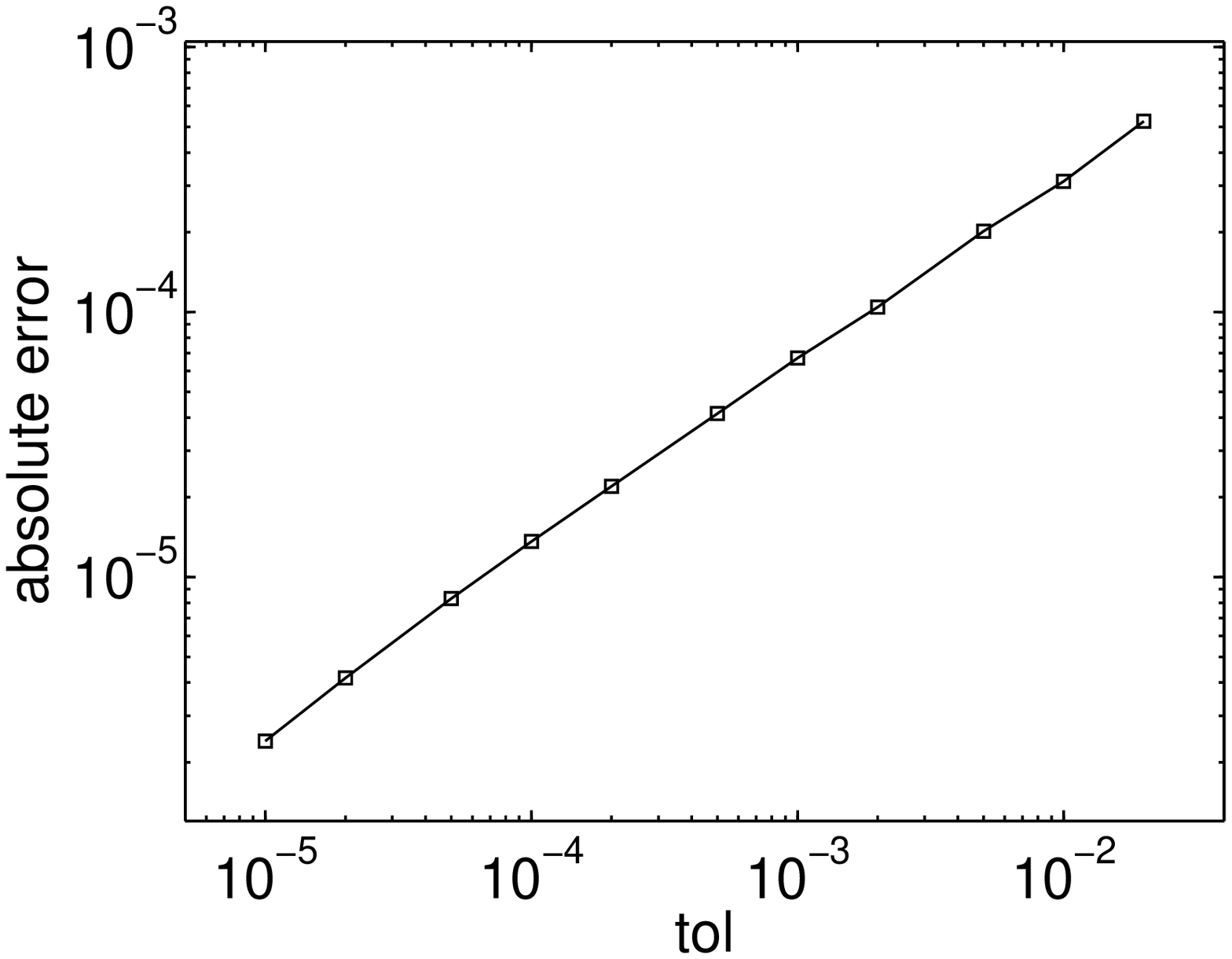}
  \\
  \includegraphics[width=0.43\textwidth, height=0.22\textheight]{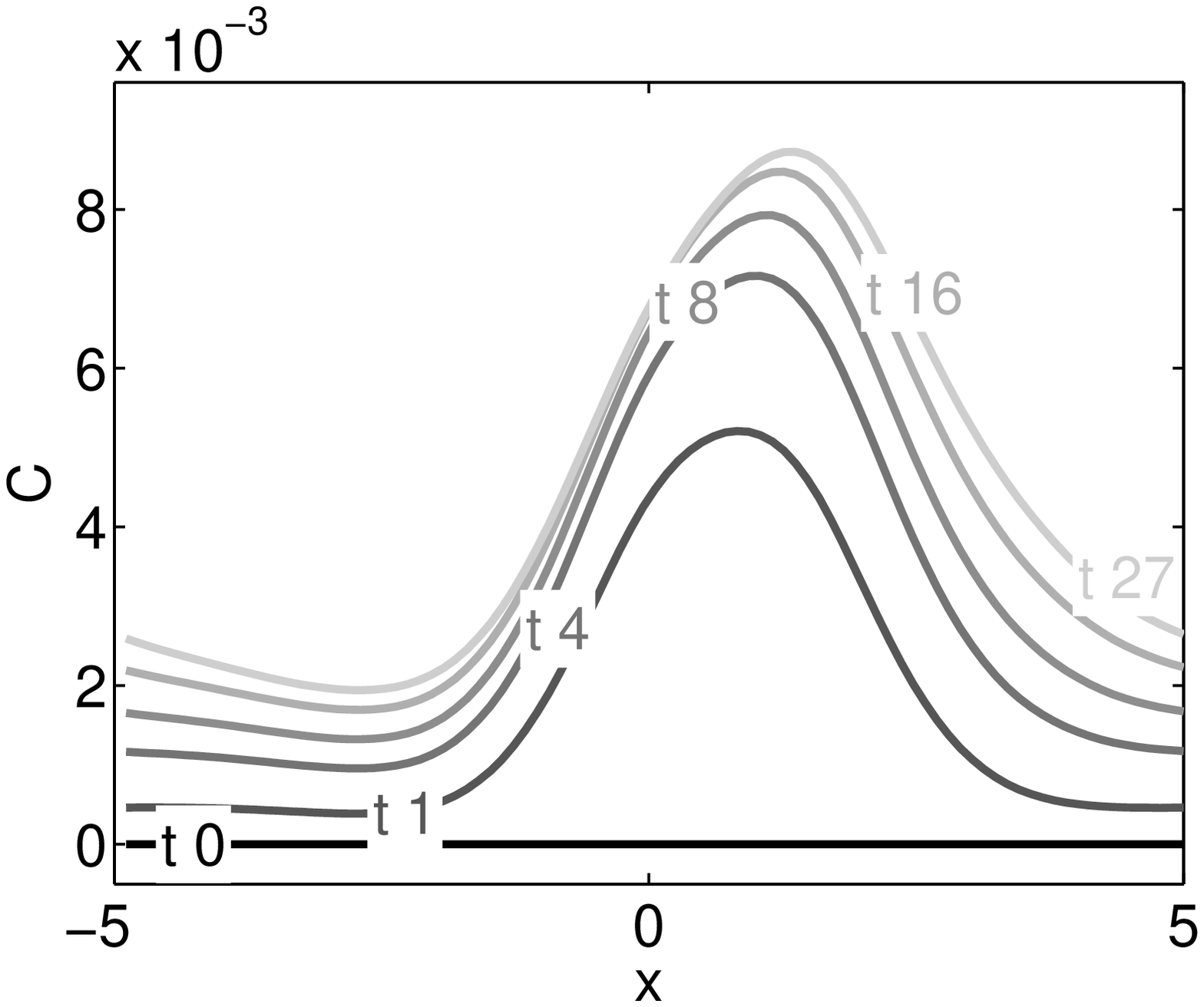}
  \includegraphics[width=0.43\textwidth, height=0.21\textheight]{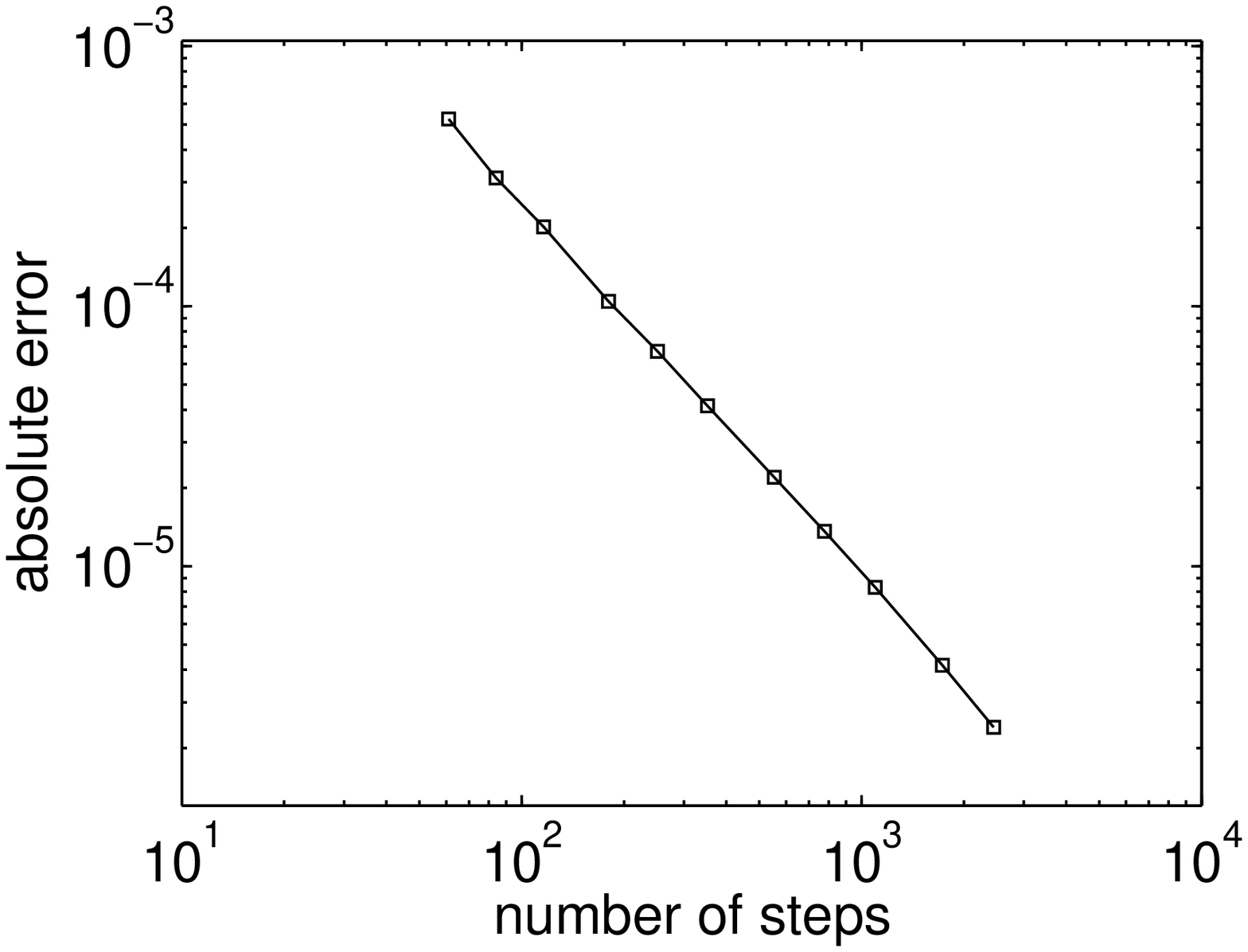}
  \caption{Left: Solutions for the three species $A$, $B$, $C$, at different times
  (lighter lines for larger times). Right: Step size versus time $t$,
  error versus tolerance, and error versus number of steps.
  All for $T=30$, $K=0.5$, $k_{1}=1$, $k_{2}=2$, and $k_{3}=3$.}
  \label{fig:err_redif_T30}
\end{figure}
Fig.~\ref{fig:err_redif_T30} shows numerical results and performance
characteristics of the algorithm. Here we inverted the Laplace
transform taking $a=1$, $d=0.5$, and $K=40$ quadrature nodes on the
hyperbolas, giving $C_1=6.036$ and $C_2=0.0739$. Again, less
stringent accuracy requirements demand fewer quadrature nodes.

\section{Dynamic fractional order viscoelasticity}

\subsection{Model}
The fractional order linear viscoelastic constitutive equation for the
stress $\Tensor{\sigma}$ considered in~\cite{AdoEL04,AdoELR06}
reads
\begin{equation}\label{eq:const_stress}
  \Tensor{\sigma} = \Tensor{\sigma}_0(t) - \gamma \int_0^t f(t- \tau)
  \Tensor{\sigma}_0 (\tau)\, d\tau,
\end{equation}
with stress tensor $\Tensor{\sigma}_{0}$ and strain tensor $\Tensor{\epsilon}$  given by
\begin{equation}
  \label{eq:stressstrain}
  \Tensor{\sigma}_{0}(t) = 2 \mu \Tensor{\epsilon}(t) + \lambda
  \mbox{tr}( \Tensor{\epsilon}(t)) \Tensor{I}
  \,; \quad
  \Tensor{\epsilon}(t) =
  \frac{1}{2} \left( \nabla \Field{u} + (\nabla \Field{u})^{T}
  \right),
\end{equation}
where $\mu$ and $\lambda$ are the Lam\'e constants and $0<
\gamma < 1$ is a given parameter. We refer to \cite{AdoEL04}
for more details about the model.

The basic equations for the displacement field $\Field{u}$ are
\begin{gather}
  \label{eq:visco}
  \begin{aligned}
    &\rho \ddot{\Field{u}}(\Svec{x},t) - \nabla \cdot
    \Tensor{\sigma}_{0}(\Field{u};\Svec{x},t)
    \\
    &
    + \gamma \int_{0}^{t} f(t-\tau)  \nabla \cdot
    \Tensor{\sigma}_{0}(\Field{u};\Svec{x},\tau) \; d\tau
    = 0 \quad \mbox{for } \Svec{x} \in \Omega; \ t \ge 0
    \\
    &\Field{u}(\Svec{x},0) =   \Field{u}_{0}(\Svec{x}) \quad \mbox{for } \Svec{x} \in
    \Omega
    \\
    &
    \dot{\Field{u}}(\Svec{x},0) =   \Field{v}_{0}(\Svec{x}) \quad \mbox{for } \Svec{x} \in
    \Omega
    \\
    &\Field{u}(\Svec{x},t) = 0 \quad \mbox{for  } \Svec{x} \in \Gamma_{D}; \ t\ge 0
    \\
    &\Tensor{\sigma}_{0}(\Field{u};\Svec{x},t) \cdot \Field{n}(\Svec{x}) =
    \Field{b}(\Svec{x},t) \quad \mbox{for  }\ \Svec{x} \in \Gamma_{N}; \ t\ge 0].
  \end{aligned}
\end{gather}
$\Gamma_{D}$ denotes the Dirichlet boundary
of the $\Omega$ and  $\Gamma_{N}$ the Neumann boundary,
where the boundary force $\Field{b}$ is applied.

On the Sobolev space
$V:=\{\Field{v} \in H^{1}(\Omega)^{2}: \Field{v}=\Field{0} \mbox{ on } \Gamma_{D} \}$
the variational formulation reads as follows:
Find $\Field{u}(t) \in V$ such that
\begin{gather}
  \label{eq:viscow}
  \begin{aligned}
    &\int_{\Omega} \rho \ddot{\Field{u}}(\Svec{x},t)\cdot \Field{\psi}(\Svec{x})
    d\Svec{x} +\int_{\Omega}\Tensor{\epsilon}(\Field{u};\Svec{x},t) : C
    \Tensor{\epsilon}(\Field{\psi}(\Svec{x})) d\Svec{x} -\int_{\Gamma_{N}}
    \Field{b}(\Svec{x},t) \cdot \Field{\psi}(\Svec{x}) d\sigma(\Svec{x})
    \\ & -
    \gamma \int_{0}^{t} f(t-\tau) \left(
    \int_{\Omega}\Tensor{\epsilon}(\Field{u};\Svec{x},\tau) : C
    \Tensor{\epsilon}(\Field{\psi}(\Svec{x})) d\Svec{x} -\int_{\Gamma_{N}}
    \Field{b}(\Svec{x},\tau) \cdot \Field{\psi}(\Svec{x}) d\sigma(\Svec{x}) \right)
    \; d\tau
    \\ &
    = 0
    \, , \quad \quad \forall \Field{\psi} \in V
    \\ &\Field{u}(0) = \Field{u}_{0} \ ;\quad
    \dot{\Field{u}}(0) = \Field{v}_{0}(\Svec{x}),
  \end{aligned}
\end{gather}
where the tensor product is given by
$$
\Tensor{\epsilon}(\Field{u}) : C \Tensor{\epsilon}(\Field{\psi})
= \sum_{i,j=1}^{2} 2\mu_{0}\epsilon_{ij}(\Field{u})\epsilon_{ij}(\Field{\psi}) + \lambda_{0}
\epsilon_{jj}(\Field{u})\epsilon_{ii}(\Field{\psi})~.
$$
Equation~\eqref{eq:viscow} is discretized in space using linear finite
elements. The mesh is generated using Triangle~\cite{Triangle} and
the assembly of the mass and stiffness matrices $M$ and $A$
and the boundary force vector $b$ is done following~\cite{AlCasterFunKlo02}.
In contrast to~\cite{AlCasterFunKlo02} we have chosen not to use Lagrange
multipliers to enforce the Dirichlet data, but to incorporate the Dirichlet data
directly. Thus~\eqref{eq:viscow} results in the abstract integro-differential equation
\begin{eqnarray*}
M \ddot{u}(t) + A u(t) - b(t) &=& \gamma \int_{0}^{t} f(t-\tau)(Au(\tau)-b(\tau))
\\
u(0) = u_{0}\ ; &&\quad \dot{u}(0) = v_{0}.
\end{eqnarray*}
The kernel $f$ in \eqref{eq:const_stress} is given by
\begin{equation}\label{kernel:visco}
f(t)= -\frac{d}{dt} E_{\alpha}\left(-t^{\alpha}\right),
\qquad 0< \alpha <1,
\end{equation}
where $E_{\alpha}$ denotes the Mittag-Leffler function of order $\alpha$, defined by
$$
E_{\alpha}(x) = \sum_{j=0}^\infty \frac{x^j}{\Gamma(1+\alpha j)}.
$$
The Laplace transform $F$ of $f$ is given by
$$
F(s)=\frac{1}{1+s^\alpha}.
$$

\subsection{Adaptive step size control}
\label{Sec:intcontroller}
The discretization of the fractional order viscoelastic problem yields
a Volterra integro-differential equation of second order of convolution type,
$$
M \ddot{u}(t) + A u(t) = \gamma\int_{0}^{t} f(t-\tau) (A u(\tau-b(\tau)) \, d\tau +
b(t) =: c(t).
$$
This is equivalent to
\begin{equation}\label{sysmatrix}
\left(
\begin{array}{c}
 \dot u\\
 \dot v
\end{array}
\right) = \left(
\begin{array}{cc}
  0 &M^{-1} \\
  -A & 0
\end{array}
\right)
 \left(
\begin{array}{c}
  v\\
  u
\end{array}
\right) +
 \left(
\begin{array}{c}
  0\\
  c
\end{array}
\right).
\end{equation}
Applying the transformations $u \to \hat u = M^{1/2}u$, $v \to \hat
v = M^{-1/2}v$, $A \to \hat A = M^{-1/2}AM^{-1/2}$ and $c \to \hat c
= M^{-1/2}c$, we get
$$
\left(
\begin{array}{c}
 \dot {\hat u}\\
 \dot {\hat v}
\end{array}
\right) = \left(
\begin{array}{c}
   \hat v\\
  -A \hat u + \hat c
\end{array}
\right).
$$
In what follows we drop the $\hat{ }s$. The time discretization is
done using the St\"ormer--Verlet scheme, which is explicit and
symmetric and has good properties for the part $\ddot{u}(t) = -A
u(t)+b(t)$ (without the memory term). The Verlet scheme for the
above equation reads
\begin{gather}
  \label{eq.Verlet}
  \begin{aligned}
    v_{n+1/2} &= v_{n} + \frac{h}{2} (-Au_{n} + c_{n})
    \\
    u_{n+1} &= u_{n} + h v_{n+1/2}
    \\
    v_{n+1} &= v_{n+1/2} + \frac{h}{2}(-Au_{n+1} + c_{n+1} ),
  \end{aligned}
\end{gather}
where $c_{n} \approx  c(t_n)$ is computed using the adaptive
convolution algorithm explained in Section~\ref{Sec:alg}. Note that
$u_{n}$ is already known before we evaluate the $c_{n}$ and thus
scheme is explicit. In order not to lose the good properties of the
St\"ormer-Verlet scheme a special step-size control is used,
following \cite{HaiLuWan06,HaiSod05} . For the \emph{integrating}
controller we fix an accuracy parameter $\veps$ (which can roughly
be viewed as the square root of a local error tolerance). Our
step-size density function should control $\ddot v$ and $\ddot u$,
therefore we take
\begin{eqnarray}
  \sigma(u,v,t) &=& \tilde{\sigma}(u,v,t)^{-1/4}
  = \bigl(\|\ddot v\|^{2} + \|\ddot u\|_{A}^{2}\bigr)^{-1/4}
   \\
   &=& \bigl((-Av+\dot c)^{T}(-Av+\dot c) +
   (-Au+c)^{T}A(-Au+c)\bigr)^{-1/4}. \nonumber
  \label{eq:sigma}
\end{eqnarray}
Assuming that $A$ is symmetric, the partial derivatives of $\tilde\sigma$ are
\begin{eqnarray*}
  \tilde\sigma_u(u,v,t) &=& 2 (Au-c)^T  AA
  \\
  \tilde\sigma_v (u,v,t) &=& 2 (Av-\dot c )^T  A
  \\
  \tilde\sigma_t (u,v,t) &=& -2(Av-\dot c)^{T} \ddot c - 2(Au-c)^{T}A\dot c.
\end{eqnarray*}
With this choice the step-size becomes approximately
proportional to $1/\sqrt{\|\ddot
u\| + c\|\ddot v\|}$. We have to take the $A$ norm of $u$ so that
$\|\ddot v\|^{2}$ and $\|\ddot u \|_{A}^{2}$ are in the same units.
We use the integrating controller of~\cite{HaiSod05},\cite[(VIII.3.2)]{HaiLuWan06}
and set
\begin{equation}
  G(u,v,t) = -\frac{1}{\sigma(u,v,t)} \nabla\sigma(u,v,t)^{T}
  \begin{pmatrix} v \\ -A u + c(t) \\ 1 \end{pmatrix}
  = -\frac{1}{4\tilde{\sigma}}\, (2(Av-\dot c)^{T} \ddot c ),
  \label{eq:Ghat}
\end{equation}
where for an evaluation, $\dot c$ and $\ddot c$ at $t_n$ are
approximated by divided differences using $c_n, c_{n-1}, c_{n-2}$
(set to zero for negative subscripts).

Transforming back to ``non-hat'' quantities and again assuming that $M$
and $A$ are symmetric one obtains
\begin{equation}
  \tilde \sigma(u,v,t)
=  (AM^{-1}v - \dot c)^{T} M^{-1} (AM^{-1}v - \dot c) +
  (Au-c)^{T}M^{-1}AM^{-1}(A u-c)
\label{eq:sigmatil}
\end{equation}
and
\begin{equation}\label{eq:Gtil}
G =
    -\frac{1}{4 \tilde{\sigma}}\,(2(AM^{-1}v-\dot c)^{T} M^{-1} \ddot c
).
\end{equation}
With an accuracy parameter $\veps$, and starting with $c_{-1}=c_{-2}=c_{0} =
b_{0}$, $z_{-1/2} = 1/\sigma(u_0,v_0,t_0) -\veps G(u_0,v_0,t_0)/2$,
we compute for $n=0,\dots $
\begin{gather}
  \label{CVerlet}
  \begin{aligned}
    z_{n+1/2} &= z_{n-1/2} + \varepsilon G(u_{n},v_{n},t_{n})
    \\
    h_{n+1/2} &= \frac{\varepsilon}{z_{n+1/2}}
    \\
    \quad t_{n+1} &= t_{n} +  h_{n+1/2}
    \\
    v_{n+1/2} &= v_{n} + \frac{h_{n+1/2}}{2}   (-Au_{n} + c_{n})
    \\
    u_{n+1} &= u_{n} + h_{n+1/2} M^{-1} v_{n+1/2}
    \\
    c_{n+1} &= \gamma \left( \frac{f_2(h_{n+1/2})}{h_{n+1/2}} (A u_{n+1}-b_{n+1}) +
    \left(f_1(h_{n+1/2})-\frac{f_2(h_{n+1/2})}{h_{n+1/2}}\right) (A u_{n} -b_{n}) \right.
    \\ & \left.+ \int_{0}^{t_{n}} f(t_{n+1}-\tau)(Au(\tau)+b(\tau)) d\tau \right) + b_{n+1}
    \\
    v_{n+1} &= v_{n+1/2} +\frac{h_{n+1/2}}{2} (-Au_{n+1} + c_{n+1}),
  \end{aligned}
\end{gather}
where $h_{n+1/2}$ is the new step-size proposed by the integrating
controller.

\subsection{Numerical example}
\begin{figure}[htbp]
  \centering
  \includegraphics[width=0.65\textwidth]{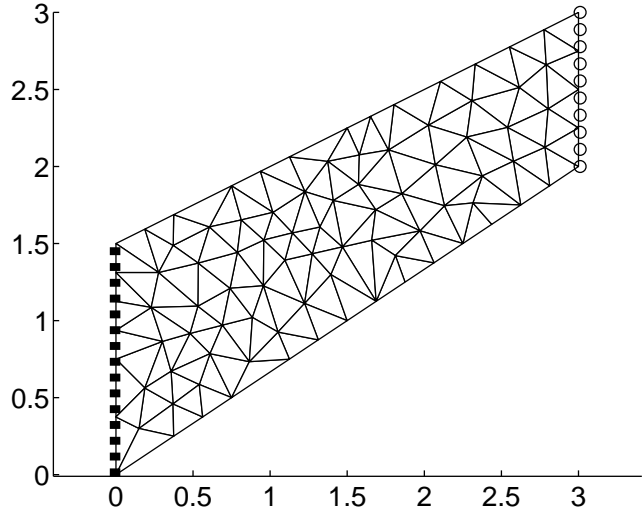}
  \caption{Cantilever with finite element mesh. $\Gamma_{D}$ is indicated by
    small squares and $\Gamma_{N_{1}}$ by circles.}
  \label{fig:viscodomain}
\end{figure}

In the example the domain $\Omega$ has the form of a cantilever as
shown in Fig~\ref{fig:viscodomain}. The Dirichlet boundary
$\Gamma_{D}$ -- the left vertical boundary of $\Omega$  -- is indicated
by squares. The time-dependent boundary force $\Field{b}$ is applied
to the right vertical boundary $\Gamma_{N_{1}}$ of $\Omega$ --
indicated by circles in Fig.~\ref{fig:viscodomain}. On
$\Gamma_{N_{2}}$ (the upper and lower part of the boundary of
$\Omega$) homogenous Neumann boundary condition is assumed.

As the initial data $\Field{u}(\Svec{x},0)$, $\Field{v}(\Svec{x},0)$ we take
the lowest mode of the static semi-discrete problem
corresponding to the first order equation \eqref{sysmatrix}, i.e., the
eigenvector of the matrix there corresponding to the smallest eigenvalue. The
boundary force is given by
$$
\Field{b}(\Svec{x},t) =
\left\{
  \begin{array}{ll}
    20 e^{1/((2t-5)^8-1)} (1,1)^{T} \quad & \mbox{for} \ 2 < t < 3;
    \Svec{x} \in \Gamma_{N_{1}}
    \\
    0  \quad & \mbox{else}.
  \end{array}
\right.
$$
In the numerical example the order of the Mittag-Leffler function is
$\alpha=1/2$. We set the density $\rho=1$ and $\gamma = 0.3$. Youngs
modulus and Poisson ratio are $E=200$, $\nu = 0.3$. Equivalently the
Lam\'e constants are $\mu = 76.9 $, $\lambda = 115.4$. For the
numerical inversion of the Laplace transform we took here $a=0.8$,
$d=0.7$ and $K=35$ quadrature nodes, giving $C_1=6.225$ and
$C_2=0.097$.

\begin{figure}
  \centering
  \includegraphics[width=0.46\textwidth,height=0.238\textheight]{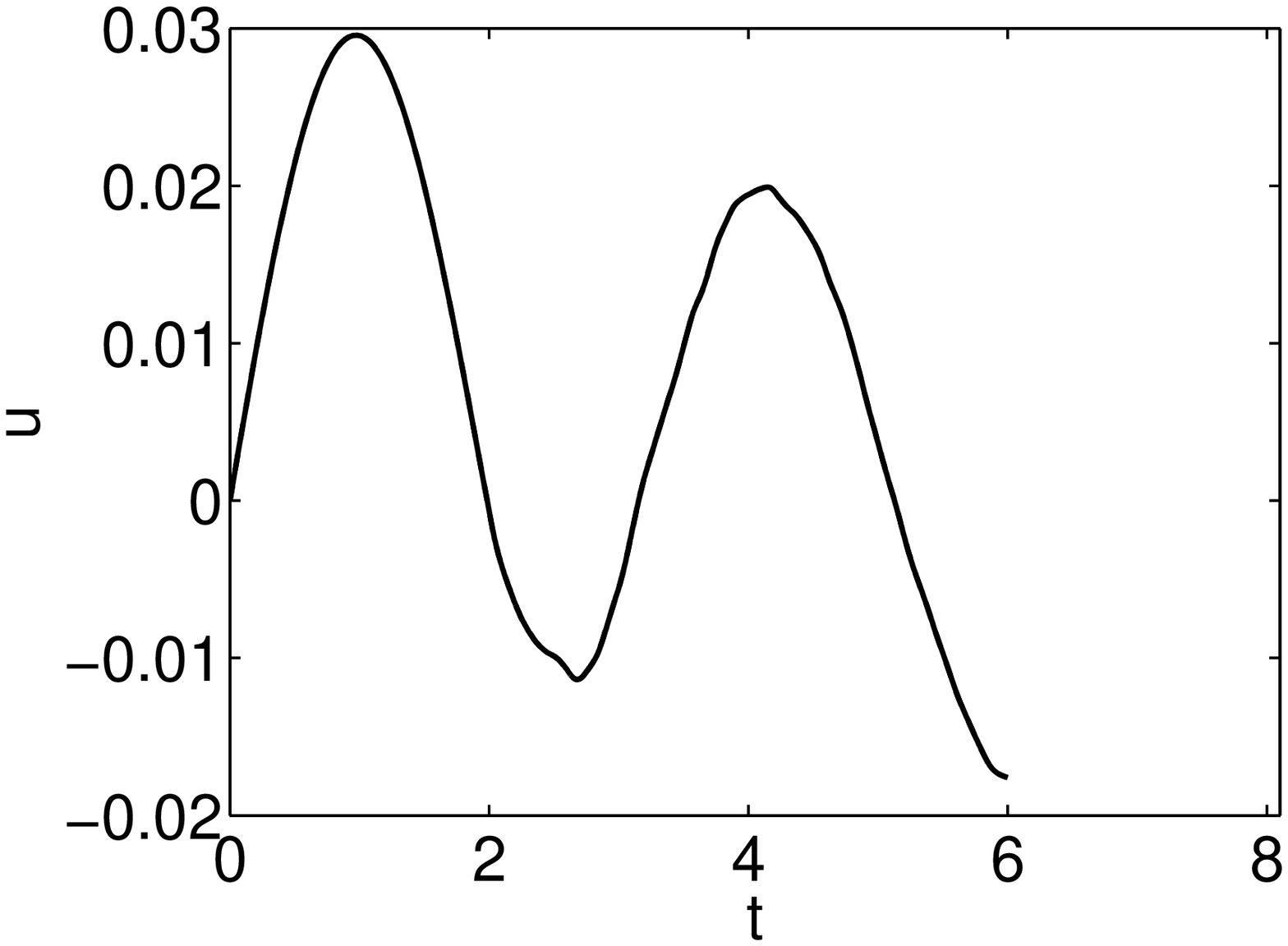}
  \includegraphics[width=0.46\textwidth,height=0.252\textheight]{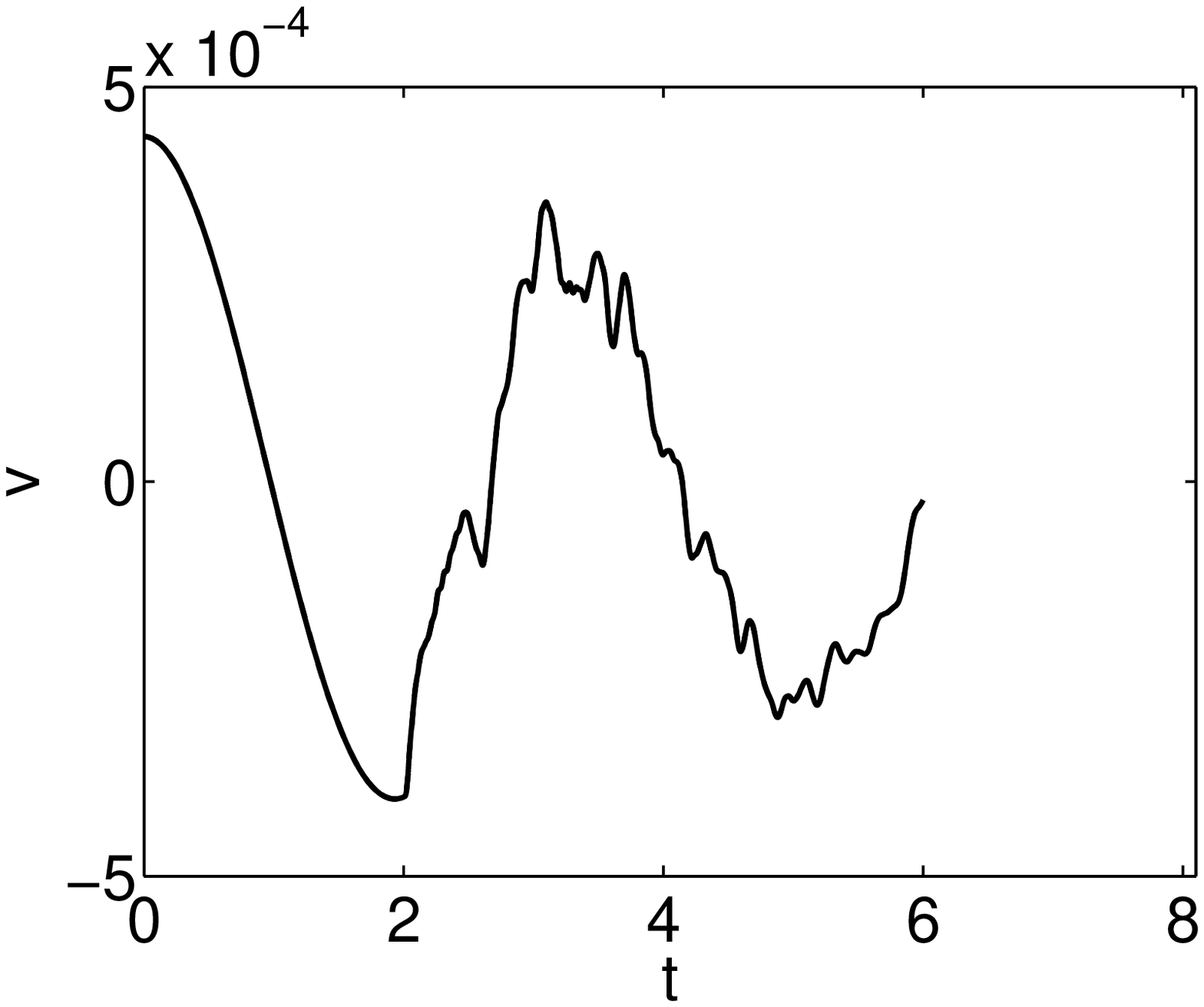}
  \\
  \includegraphics[width=0.46\textwidth,height=0.253\textheight]{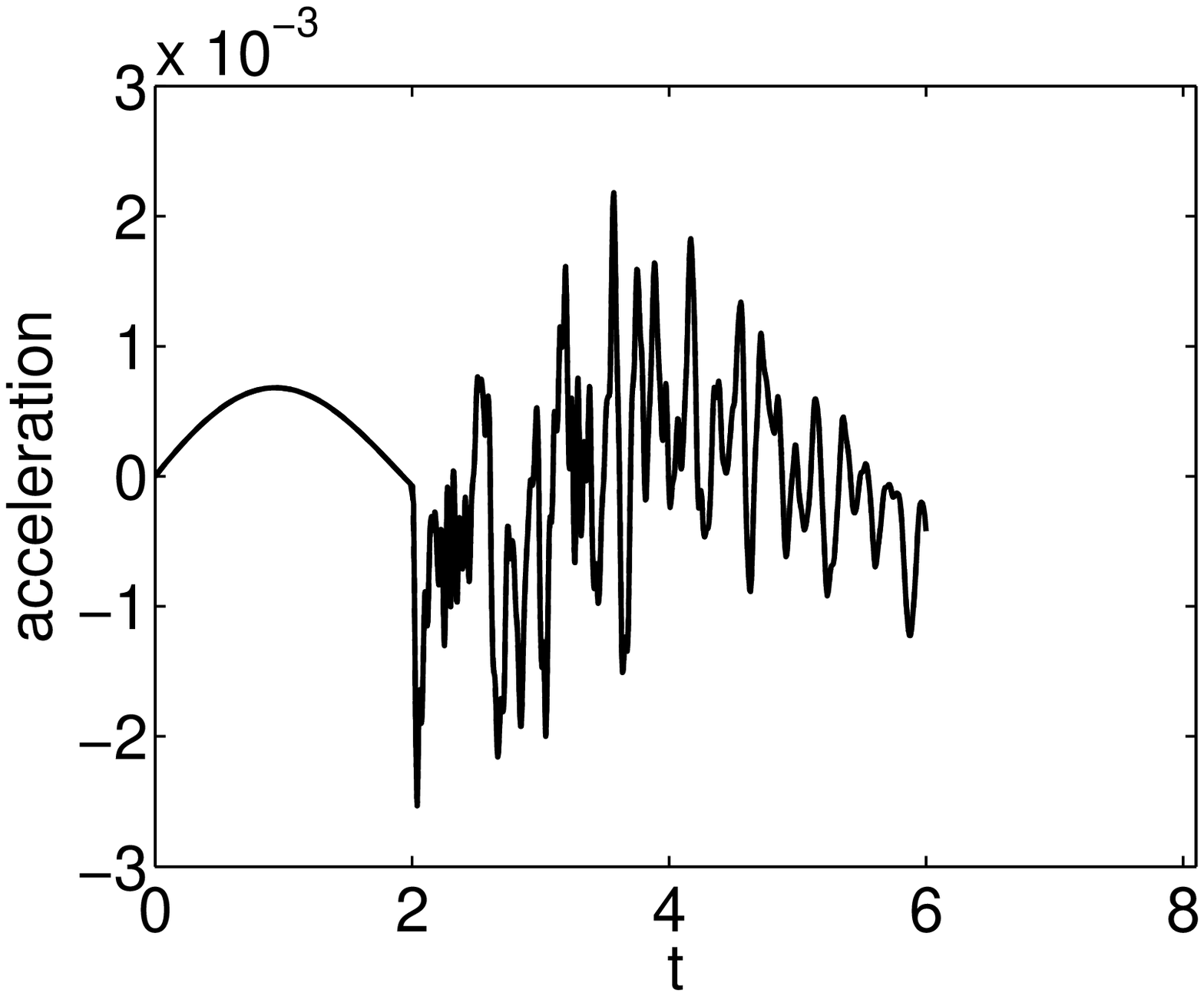}
  \includegraphics[width=0.46\textwidth,height=0.231\textheight]{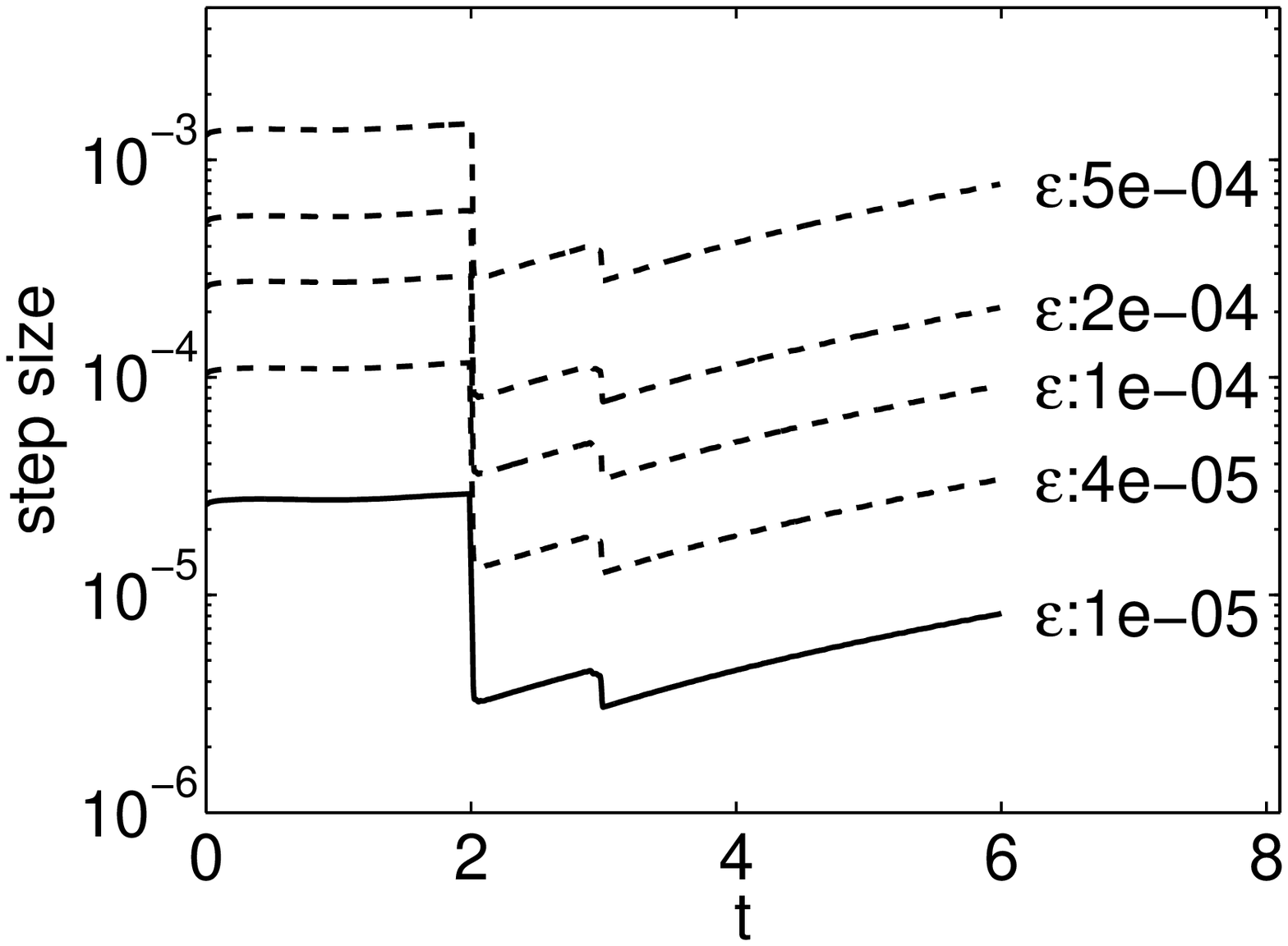}
  \caption{$t=6$. Displacement, velocity, and acceleration recorded for one
  fixed spatial node. Also step size versus time for different $\veps$.}
  \label{fig:evolelasto}
\end{figure}

The evolution from $t=0$ to $t=6$ of the horizontal component ($x$ coordinate)
of the displacement field $\Field{u}$, the velocity $\dot{\Field{u}}$, and the
acceleration $\ddot{\Field{u}}$ recorded at the upper left corner of the
cantilever is shown in Fig.~\ref{fig:evolelasto}.
At $t=2$ when the boundary force is applied, an abrupt
change in the velocity and strong oscillations in the acceleration is observed.
Furthermore Fig.~\ref{fig:evolelasto} shows the evolution
of the step-size $h_{n+1/2}$ for five different precision parameters $\veps
= 1\cdot10^{-5}, 4\cdot 10^{-5}, 1\cdot 10^{-4}, 2\cdot 10^{-4}, 5\cdot 10^{-4}$.
The integrating controller reduces the step-size by roughly a factor of ten at
$t=2$.

\begin{figure}
  \centering
  \includegraphics[width=0.46\textwidth]{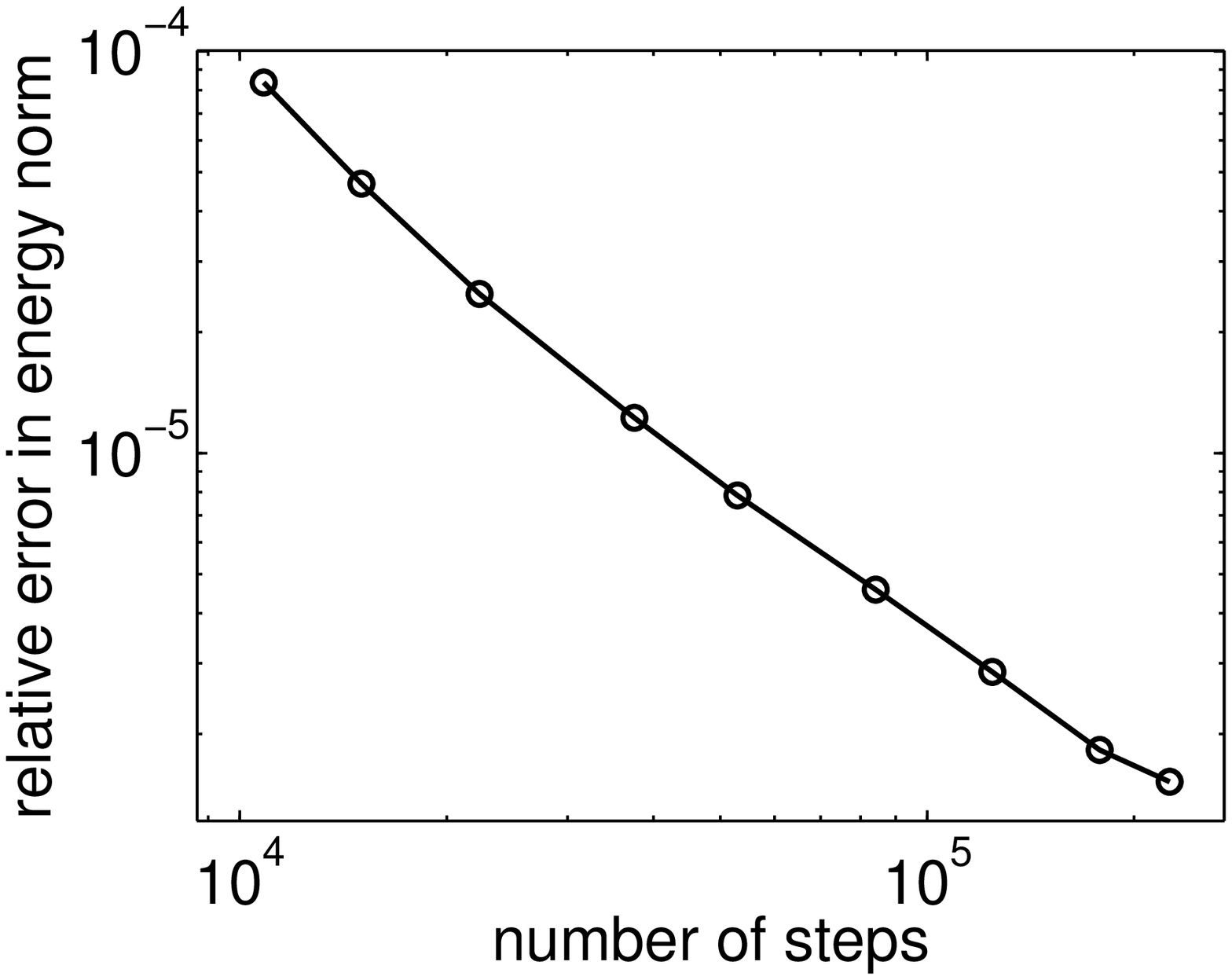}
  \includegraphics[width=0.46\textwidth]{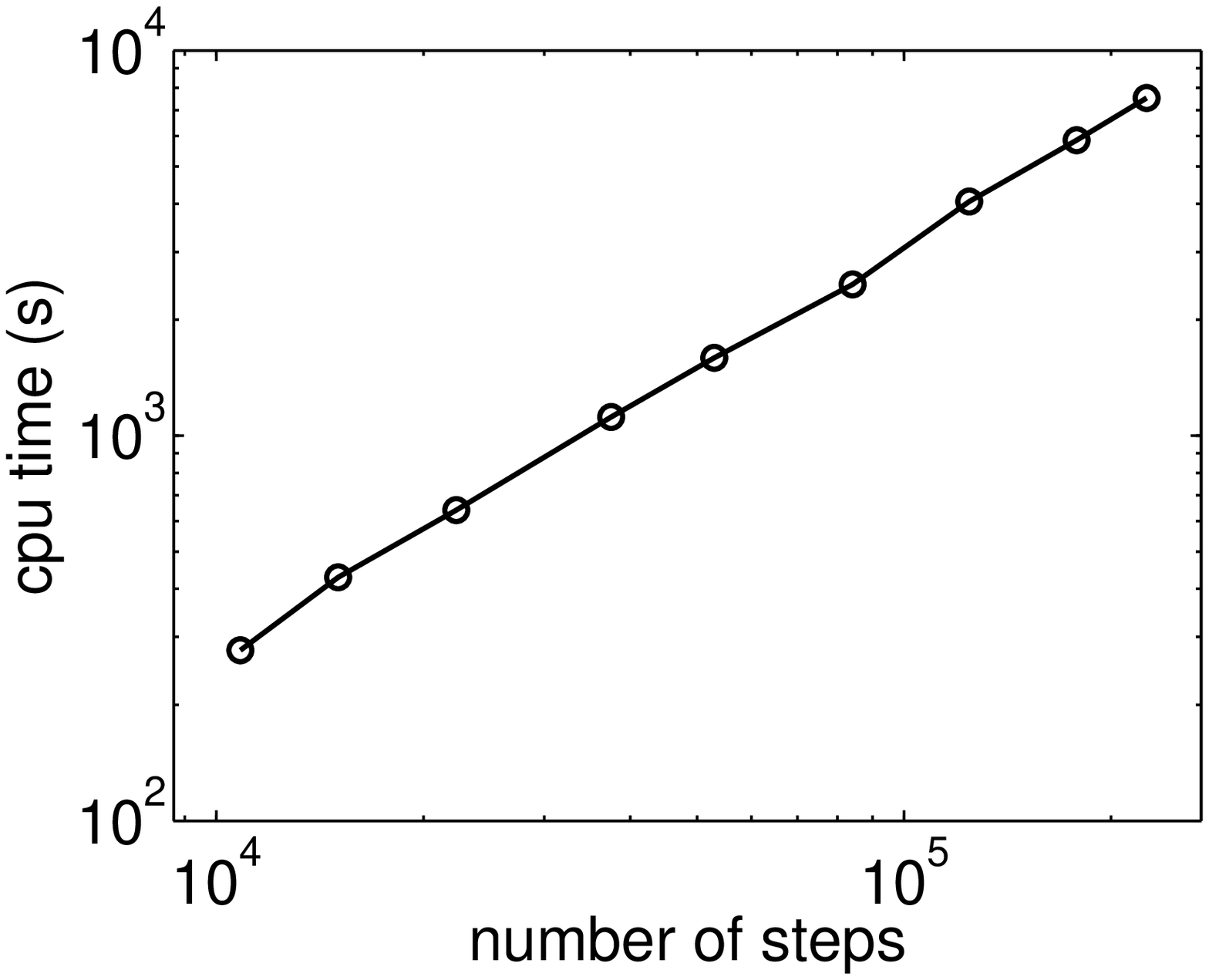}
  \caption{Error versus number of steps and error versus cpu time, for $t=6$.}
  \label{fig:errorelasto}
\end{figure}

The relative error at $t=6$ in energy norm $\|u \|_{A} + \|v\|_{M}$
is calculated with respect to a reference solution obtained with
$\veps = 1\cdot10^{-5}$. The left plot in Fig~\ref{fig:errorelasto}
shows the relative error versus the number of steps. The different
solutions were obtained for $\veps = 4\cdot 10^{-5}, 5\cdot 10^{-5},
7\cdot 10^{-5}, 1\cdot 10^{-4}, 1.5\cdot 10^{-4}, 2\cdot 10^{-4},
3\cdot 10^{-4}, 4\cdot 10^{-4}, 5\cdot 10^{-4}$. The right plot in
Fig~\ref{fig:errorelasto} shows the cpu time in seconds versus the
number of steps, illustrating the essentially linear growth.

\section*{Acknowledgment}
We thank Roland Klose for kindly providing his elasticity finite-element code.
\bibliographystyle{amsplain}

\appendix

\section{Pseudocodes for the algorithm}
We describe one step of the algorithm from $t_{n-1}$ to a given new time $t=t_{n}$.

For all $\ell=1,\dots,L$ and $k=-K,\dots,K$ the ODEs~\eqref{eq:ode}
corresponding to the $\lambda_k^{(\ell)}$ with initial time
$t_{\ell}^-$ are advanced to the new time $t$ or restarted,
depending on whether the horizontal line at height $t$ fits the
current patch or enters a new patch (see
Section~\ref{Subsec:filling} and Figure \ref{fig:mosaicYM}), i.e.,
we compute $y^{(\ell)}(t_{\ell}^-,t,\lambda_k^{(\ell)})$, or we set
$t_{\ell}^-=t$ and $y^{(\ell)}(t_{\ell}^-,t,\lambda_k^{(\ell)}) =
0$, using the pseudocode~\ref{odesol}.

In the following pseudocodes \texttt{Y(l)} denotes a structure
storing:
\begin{itemize}
\item $\texttt{Y(l).data} \leftarrow y$, the solution of
  the ODE corresponding to the $\ell$-patch,
\item $\texttt{Y(l).tini} \leftarrow t_{\ell}^-$, the initial time of the ODE in the $\ell$-patch,
\item $\texttt{Y(l).tcur} \leftarrow t_{n}$, the new time (new final time of the ODE),
\item $\texttt{Y(l).b} \leftarrow b_{\ell}$, the number of the current step
  in the corresponding $\ell$-patch of the mosaic, along the vertical
  line from $(t_n,0)$ to $(t_n,t_n)$ (see Figure~\ref{fig:mosaicYM}).
  This value ranges between 1 and $B$. In the example of
  Section~\ref{Subsec:filling} we have for $t_n=t_{15}$, $b_3=3$,
  $b_2=1$ and $b_1=3$.
\item $\texttt{Y(l).tmin} \leftarrow t_{min}^{(\ell)} =
  h_{min} \sum_{k=\ell+1}^{L} b_{k} B^{k-1}\,$, the baseline of the
  $\ell$-patch of the mosaic, along the vertical line from $(t_n,0)$
  to $(t_n,t_n)$. 
\item $\texttt{Y(l).tmax} \leftarrow t_{max}^{(\ell)} =
  h_{min} \sum_{k=\ell}^{L} b_{k} B^{k-1} \,$, the top line of the
  $\ell$-patch of the mosaic, along the vertical line from $(t_n,0)$
  to $(t_n,t_n)$ 
  ($t_{max}^{(\ell+1)}=t_{min}^{(\ell)}$),
\item  $\texttt{Y(l).ub} \leftarrow  h_{min} \left(1+ \sum_{k=0}^{\ell} B^k\right)$
  upper bound of the approximation interval,
\item  $\texttt{Y(l).lb} \leftarrow  h_{min} \left(1+ \sum_{k=0}^{\ell-2} B^k\right)$
  lower bound of the approximation interval,
\item $\texttt{Y(l).gini} \leftarrow g(t_{\ell}^-)$ and $\texttt{Y(l).gcur}
  \leftarrow g(t_{n})$, the values of the in\-ho\-mo\-ge\-nei\-ty at
  the initial time in the $\ell$-patch and at the current time.
\end{itemize}
In the rare case that $t_{n}$ is exactly on one of the horizontal
lines of the mosaic, the structure \texttt{Y(l)} is copied to
\texttt{YA(l)} before restarting for bookkeeping purposes.
\begin{algorithm}[H]
  \caption{Advance and restart the scalar ODEs. \texttt{odesol}}
  \label{odesol}
  \begin{algorithmic}
    \FOR{$l = 1$ to $L$}
      \IF{$t_{n} \ge Y(l).tmin + Y(l).b*B^{(l-1)}*hmin$}
        \IF{$ t_{n} \ge Y(l).tmax$}
          \IF{$t_{n} = Y.tmax$}
          \STATE $Y(l) =
          \mathtt{odesadvance}(Y(l),t_{n},g_{n-1},g_{n})$ {\bf ;}
          \STATE $ Y(l).b = B $ {\bf ;} \STATE $ YA(l) = Y(l)$ {\bf ;}
          \ENDIF
          \STATE restart the ODE $Y(l)$ \textit{c.f.}
          Algorithm~\ref{restart} {\bf ;}
        \ELSE
          \STATE $Y(l) = \mathtt{odesadvance}(Y(l),t_{n},g_{n-1},g_{n})$ {\bf ;}
          \WHILE{$t > Y(l).tmin+Y(l).b*B^{(l-1)}*hmin$}
            \STATE $Y(l).b = Y(l).b + 1$ {\bf ;}
          \ENDWHILE
        \ENDIF
      \ELSE
        \STATE $Y(l) =\mathtt{odesadvance}(Y(l),t_{n},g_{n-1},g_{n})$ {\bf ;}
      \ENDIF
    \ENDFOR
    \end{algorithmic}
\end{algorithm}
\begin{algorithm}[H]
  \caption{restart the ODE}
  \label{restart}
  \begin{algorithmic}
  \STATE $Y(l).data = 0*Y(l).data$ {\bf ;} \quad $ Y(l).b = 1$ {\bf ;}
  \WHILE{$t_{n} > Y(l).tmin+Y(l).b*B^{(l-1)}*hmin$}
    \STATE $ Y(l).b = Y(l).b + 1$ {\bf ;}
  \ENDWHILE
  \STATE $Y(l).tini = t${\bf ;} \quad $ Y(l).tcur = t${\bf ;}
  \quad  $Y(l).gini = g_{n}${\bf ;}
  \quad $Y(l).gcur = g_{n}$ {\bf ;}
  \WHILE{$t_{n} \ge Y(l).tmax$}
    \STATE $Y(l).tmin = Y(l).tmax$ {\bf ;}
    \STATE $ Y(l).tmax = Y(l).tmin+B^l*hmin$ {\bf ;}
  \ENDWHILE
  \end{algorithmic}
\end{algorithm}
\begin{algorithm}[H]
  \caption{Advancing the ODE (exponential Euler method) \texttt{odesadvance}}
  \label{odeadvance}
  \begin{algorithmic}
    \FOR{$k=-K,\dots,K$}
    \STATE $Y(\ell).data_{k} = Y(\ell).data_{k} +
    (\exp(dt*\lambda_{k}^{(\ell)})-1)/\lambda_{k}^{(\ell)}$
    \STATE \hspace*{3em} $*(\lambda_{k}^{(\ell)}*Y(\ell).data_{k}+g_{n-1}+(g_{n}-g_{n-1})/dt/\lambda_{k}^{(\ell)})
    + (g_{n-1}-g_{n})/\lambda_{k}^{(\ell)}$ {\bf ;}
    \ENDFOR
    \STATE $Y(\ell).gcur = g_{n}$ {\bf ;} \quad $Y(\ell).tcur = t$ {\bf ;}
  \end{algorithmic}
\end{algorithm}

To fill the mosaic botton-up, Algorithm~\ref{update} is used. There copying
the structure  $\texttt{Y(l)}$ to
$\texttt{YM(l)}$ and $\texttt{YA(l)}$ is done
by checking if the distances to the diagonal
$t_n-\texttt{Y(l)}.t_{min}$ and $t_n-\texttt{Y(l)}.t_{max}$ fit the
approximation interval $I_{\ell}$.
\begin{algorithm}[H]
  \caption{update routine \texttt{update}}
  \label{update}
  \begin{algorithmic}
    \FOR{$\ell$ from $L$ downto $1$}
    \IF{$t_{n}-YM(\ell).tmin-hmin*YM(\ell).b*B^{(\ell-1)} \ge Y(\ell).lb$ \&
      $t_{n}-YM.tini \le Y(\ell).ub$}
      \STATE $YT(\ell) = YM(\ell)$ {\bf ;}
    \ENDIF
    \IF{$t_{n} \ge Y(\ell).tmax$}
      \STATE $YA(\ell) = Y(\ell)$ {\bf ;} \STATE  $YM(\ell) = Y(\ell)$ {\bf ;}
    \ELSIF{$t_{n} \ge Y(\ell).tmin + hmin*YM.b*B^{(\ell-1)}$}
      \STATE $YM(\ell) = Y(\ell)$ {\bf ;}
    \ENDIF
    \IF{$t_{n}-YM(\ell).tmin-hmin*YM(\ell).b*B^{(\ell-1)} \ge Y(\ell).lb $ \&
      $ t_{n}-YM.tini \le Y(\ell).ub$}
      \STATE $YT(\ell) = YM(\ell)$ {\bf ;}
    \ENDIF
    \IF{$t_{n}-YA(\ell).tmin-hmin*YA(\ell).b*B^{(\ell-1)} \ge Y(\ell).lb $ \&
      $ t_{n}-YA.tini \le Y(\ell).ub $ \&
      $ YA(\ell).tini > YT(\ell).tini$}
      \STATE $YT(\ell) = YA(\ell)$ {\bf ;}
    \ENDIF
  \ENDFOR
  \end{algorithmic}
\end{algorithm}
\texttt{YT(l)} is updated by
checking if either of the distances to the diagonal corresponding to
\texttt{YM(l)} or \texttt{YA(l)} fit the approximation intervals.
\begin{algorithm}[H]
  \caption{update routine part 2 \texttt{update}}
  \label{update2}
  \begin{algorithmic}
    \STATE $idv = [\ ]$ {\bf ;} \quad $ito = 1$ {\bf ;}
   \IF{$t_{n}-YT(L).tcur \ge Y(L).lb $ \&
      $ t_{n}-YT(L).tini \le Y(L).ub$}
    \STATE $idv = [L,\ idv]$ {\bf ;}
    \ELSE
    \STATE $YT(L).tini = -inf$ {\bf ;} \quad $YT(L).tcur = -inf$ {\bf ;}
  \ENDIF
  \FOR{$\ell$ from $L-1$ downto $1$}
    \IF{$ t_{n}-YT(\ell).tcur \ge Y(\ell).lb $ \&
      $ t_{n}-YT(\ell).tini \le Y(\ell).ub $ \&
      $ YT(\ell).tini > YT(\ell+1).tini $ \&
      $ YT(\ell).tcur > YT(\ell+1).tcur$}
      \STATE $idv = [\ell,\ idv]$ {\bf ;}
   \ELSE
     \STATE $YT(\ell).tini = -inf$ {\bf ;}
     \quad  $YT(\ell).tcur = -inf$ {\bf ;}
    \ENDIF
    \IF{$t_{n}-t_{n-1} \ge Y(\ell).lb$ \&
      $t_{n}-t_{n-1} \le Y(\ell).ub$ }
     \STATE $ito = \ell$ {\bf ;}
    \ENDIF
  \ENDFOR
  \STATE $idv = [ito\ idv]$  {\bf ;}
  \end{algorithmic}
\end{algorithm}
Algorithm~\ref{fastconvquad} puts together the necesary ``direct''
and ``odes'' steps.
It uses the ode solutions \texttt{YT} and the vector \texttt{idv} from
Algorithm~\ref{update2}. \texttt{idv} stores the orders $\ell$ of the
$\texttt{YT}$ required at $t= t_{n}$.
\begin{algorithm}[H]
  \caption{Fast convolution evaluation, \textit{cf.}~\eqref{eq:convsplit} }
  \label{fastconvquad}
  \begin{algorithmic}
    \STATE $out = 0$ {\bf ;}
    \IF{$idv$ is not empty}
      \IF{$\mbox{length}(idv) \ge 2$}
        \IF{$t_{n-1} > YT(idv(2)).tcur$}
          \STATE $out = out + $
          \STATE $\mathtt{directstep}(t,YT(idv(2)).tcur, t_{n-1},YT(idv(2)).gcur,
          g_{n-1}
          )$ {\bf ;} \eqref{eq:linapprox}
        \ENDIF
        \FOR{$ll = 2:length(idv)-1$}
          \IF{$YT(idv(ll)).tcur  > YT(idv(ll)).tini$}
            \STATE $\ell = idv(ll)$ {\bf ;}
            \STATE $
            out = out + \sum_{k=-K}^{K} w^{(\ell)}_{k} F(\lambda^{(\ell)}_{k})
            \exp(t-YT(\ell).tcur)*
            YT(\ell).data_{k}$ {\bf ;} \eqref{eq:approx_convode}
          \ENDIF
          \IF{$YT(idv(ll)).tini > YT(idv(ll+1)).tcur$}
            \STATE $out = out + \mathtt{directstep}(t,
            YT(idv(ll+1)).tcur, YT(idv(ll)).tini, $
            \STATE \hspace*{3em}$YT(idv(ll+1)).gcur,YT(idv(ll)).gini
            )$ {\bf ;} \eqref{eq:linapprox}
          \ENDIF
        \ENDFOR
        \IF{$YT(idv(end)).tcur  > YT.(idv(end))tini$}
          \STATE $\ell = idv(end)$ {\bf ;}
          \STATE$
          out = out + \sum_{k=-K}^{K} w^{(\ell)}_{k} F(\lambda^{(\ell)}_{k})
          \exp(t-YT(\ell).tcur)*YT(\ell).data_{k} $ {\bf ;} \eqref{eq:approx_convode}
        \ENDIF
      \ENDIF
    \ENDIF
  \end{algorithmic}
\end{algorithm}

\end{document}